\documentclass[11pt,oneside]{amsart}

\usepackage{fullpage}

\usepackage{amsthm,amsfonts,amssymb,amsmath,amsxtra}
\usepackage[all]{xy}
\SelectTips{cm}{}
\usepackage{tikz}
\usepackage{verbatim}
\usepackage{mathrsfs}

\usepackage{mathabx}

\RequirePackage{xspace}
\RequirePackage{etoolbox}
\RequirePackage{varwidth}
\RequirePackage{enumitem}
\RequirePackage{tensor}
\RequirePackage{mathtools}
\RequirePackage{longtable}
\RequirePackage{multirow}
\RequirePackage{makecell}
\RequirePackage{bigints}

\usepackage[backref=page]{hyperref} %

\usepackage[T2A,T1]{fontenc}
\DeclareSymbolFont{cyrillic}{T2A}{cmr}{m}{n}
\DeclareMathSymbol{\Sha}{\mathalpha}{cyrillic}{216}

\newcommand{\sA}{\ensuremath{\mathscr{A}}\xspace}

\newcommand{\sS}{\ensuremath{\mathscr{S}}\xspace}

\newcommand{\Ch}{{\mathrm{CH}}}
\newcommand{\aCh}{{\widehat{\mathrm{CH}}}}

\newcommand{\disc}{{\mathrm{disc}}}

\DeclareMathOperator{\End}{End}

\newcommand{\GL}{\mathrm{GL}}

\DeclareMathOperator{\Hom}{Hom}

\let\Im\relax

\DeclareMathOperator{\Im}{Im}
\newcommand{\Ind}{{\mathrm{Ind}}}

\DeclareMathOperator{\Int}{\ensuremath{\mathrm{Int}}\xspace}

\DeclareMathOperator{\Lie}{Lie}

\newcommand{\Mp}{{\mathrm{Mp}}}

\DeclareMathOperator{\ord}{ord}

\DeclareMathOperator{\rank}{rank}

\renewcommand{\Re}{{\mathrm{Re}}}

\DeclareMathOperator{\Res}{Res}

\newcommand{\Sh}{\mathrm{Sh}}

\DeclareMathOperator{\Spec}{Spec}

\newcommand{\Sp}{{\mathrm{Sp}}}

\DeclareMathOperator{\tr}{tr}

\renewcommand{\U}{\mathrm{U}}
\renewcommand{\O}{\mathrm{O}}

\DeclareMathOperator{\vol}{vol}

\newcommand{\wit}{\widetilde}
\newcommand{\wh}{\widehat}

\newcommand{\lra}{\longrightarrow}

\setitemize[0]{leftmargin=*,itemsep=\the\smallskipamount}
\setenumerate[0]{leftmargin=*,itemsep=\the\smallskipamount}

\renewcommand{\to}{%
   \ifbool{@display}{\longrightarrow}{\rightarrow}%
   }
\let\shortmapsto\mapsto
\renewcommand{\mapsto}{%
   \ifbool{@display}{\longmapsto}{\shortmapsto}%
   }
\newlength{\olen}
\newlength{\ulen}
\newlength{\xlen}
\newcommand{\xra}[2][]{%
   \ifbool{@display}%
      {\settowidth{\olen}{$\overset{#2}{\longrightarrow}$}%
       \settowidth{\ulen}{$\underset{#1}{\longrightarrow}$}%
       \settowidth{\xlen}{$\xrightarrow[#1]{#2}$}%
       \ifdimgreater{\olen}{\xlen}%
          {\underset{#1}{\overset{#2}{\longrightarrow}}}%
          {\ifdimgreater{\ulen}{\xlen}%
             {\underset{#1}{\overset{#2}{\longrightarrow}}}
             {\xrightarrow[#1]{#2}}}}%
      {\xrightarrow[#1]{#2}}
   }
\makeatother
\newcommand{\xyra}[2][]{%
   \settowidth{\xlen}{$\xrightarrow[#1]{#2}$}%
   \ifbool{@display}%
      {\settowidth{\olen}{$\overset{#2}{\longrightarrow}$}%
       \settowidth{\ulen}{$\underset{#1}{\longrightarrow}$}%
       \ifdimgreater{\olen}{\xlen}%
          {\mathrel{\xymatrix@M=.12ex@C=3.2ex{\ar[r]^-{#2}_-{#1} &}}}%
          {\ifdimgreater{\ulen}{\xlen}%
             {\mathrel{\xymatrix@M=.12ex@C=3.2ex{\ar[r]^-{#2}_-{#1} &}}}
             {\mathrel{\xymatrix@M=.12ex@C=\the\xlen{\ar[r]^-{#2}_-{#1} &}}}}}%
      {\mathrel{\xymatrix@M=.12ex@C=\the\xlen{\ar[r]^-{#2}_-{#1} &}}}%
   }
\makeatletter
\newcommand{\xla}[2][]{%
   \ifbool{@display}%
      {\settowidth{\olen}{$\overset{#2}{\longleftarrow}$}%
       \settowidth{\ulen}{$\underset{#1}{\longleftarrow}$}%
       \settowidth{\xlen}{$\xleftarrow[#1]{#2}$}%
       \ifdimgreater{\olen}{\xlen}%
          {\underset{#1}{\overset{#2}{\longleftarrow}}}%
          {\ifdimgreater{\ulen}{\xlen}%
             {\underset{#1}{\overset{#2}{\longleftarrow}}}
             {\xleftarrow[#1]{#2}}}}%
      {\xleftarrow[#1]{#2}}
   }
\newcommand{\isoarrow}{%
   \ifbool{@display}{\overset{\sim}{\longrightarrow}}{\xrightarrow\sim}%
   }
\renewcommand{\lra}{%
   \ifbool{@display}{\longleftrightarrow}{\leftrightarrow}%
   }

\newcommand{\Herm}{\mathrm{Herm}}
\newcommand{\rd}{\mathrm{d}}

\newcommand{\KG}{K}

\newcommand{\pEis}{\partial\mathrm{Eis}}
\newcommand{\Eis}{\mathrm{Eis}}

\newcommand{\wdeg}{\widehat\deg}

\newcommand{\sx}{\mathsf{x}}
\newcommand{\sy}{\mathsf{y}}

\DeclareSymbolFont{sfletters}{OML}{cmbrm}{m}{it}
\DeclareMathSymbol{\stau}{\mathord}{sfletters}{"1C}
\newcommand{\sz}{\tau}

\DeclareFontFamily{U}{matha}{\hyphenchar\font45}
\DeclareFontShape{U}{matha}{m}{n}{
      <5> <6> <7> <8> <9> <10> gen * matha
      <10.95> matha10 <12> <14.4> <17.28> <20.74> <24.88> matha12
      }{}
\DeclareSymbolFont{matha}{U}{matha}{m}{n}
\DeclareFontFamily{U}{mathx}{\hyphenchar\font45}
\DeclareFontShape{U}{mathx}{m}{n}{
      <5> <6> <7> <8> <9> <10>
      <10.95> <12> <14.4> <17.28> <20.74> <24.88>
      mathx10
      }{}
\DeclareSymbolFont{mathx}{U}{mathx}{m}{n}

\DeclareMathSymbol{\obot}         {2}{matha}{"6B}

\newtheorem{theorem}[subsubsection]{Theorem}

\newtheorem {conjecture}[subsubsection]{Conjecture}

\theoremstyle{definition}

\newtheorem{example}[subsubsection]{Example}

\newtheorem{remark}[subsubsection]{Remark}

\numberwithin{equation}{subsubsection}
\newtheorem{assumption}[subsubsection]{Assumption}

\newcommand{\etale}{\'etale }

\newcommand{\m}{n}

\newcommand{\Blue}{\emph}
\newcommand{\DL}{\mathop {\operator@font DL}\nolimits}
\newcommand{\UU}{\mathrm{U}}

\newcommand{\Gn}{G}

\newcommand{\Mn}{M}
\newcommand{\Nn}{N}
\newcommand{\VA}{V(\mathbb{A})}
\newcommand{\VAf}{\mathbb{V}_f}
\newcommand{\n}{m}
\newcommand{\mm}{n}

\renewcommand{\H}{\mathrm{H}}
\newcommand{\Ah}{\sA_{m}(G(\mathbb{A}))}
\newcommand{\KM}{\mathrm{KM}}
\newcommand{\HX}{\H^{2\mm}(X(\mathbb{C}),\mathbb{C})}
\newcommand{\HXK}{\H^{2\mm}(X_K(\mathbb{C}),\mathbb{C})}
\newcommand{\CHX}{\Ch^{\mm}(X)_\mathbb{C}}
\newcommand{\CHXK}{\Ch^{\mm}(X_K)_\mathbb{C}}
\newcommand{\XK}{\mathcal{X}_K}
\newcommand{\voln}{\vol^\natural}

\newcommand{\varphit}{\wit\varphi^V}
\newcommand{\varphic}{\varphi^V}
\newcommand{\varphiK}{\varphi}
\newcommand{\varphiI}{\varphi^\mathbb{V}}

\newcommand{\hZ}{\widehat{\mathcal{Z}}(\sz,\varphi)}
\newcommand{\aCHXK}{\aCh^\m(\mathcal{X}_K)}
\newcommand{\avol}{\wh\vol}
\newcommand{\ceq}{\stackrel{\cdot}{=}}

\title{Geometric and arithmetic theta correspondences}

\author[Chao Li]{Chao Li}
\address{Columbia University, Department of Mathematics, 2990 Broadway,	New York, NY 10027, USA}
\email{chaoli@math.columbia.edu}

\date{February 24, 2024}


\thanks{The author would like to thank the organizers of the IHES summer school (Pierre-Henri Chaudouard, Wee Teck Gan, Tasho Kaletha and Yannis Sakellaridis)  for the invitation and their tremendous effort to make the summer school successful. The author is grateful to Wee Teck Gan for his help, and to Wei Zhang and the anonymous referee for the careful reading and many helpful comments. The author's work is partially supported by the NSF grant DMS-2101157.}

\begin{document}

\maketitle{}

\begin{abstract}
Geometric/arithmetic theta correspondences provide correspondences between automorphic forms and cohomology classes/algebraic cycles on Shimura varieties. We give an introduction focusing on the example of unitary groups and highlight recent advances in the arithmetic theory (also known as the Kudla program) and their applications.  These are expanded lecture notes
   for the IHES 2022 Summer School on the Langlands Program. 
\end{abstract}

\section{Introduction}

In the 1990s, Kudla \cite{Kudla1997a} initiated a far-reaching program relating arithmetic geometry (special cycles on Shimura varieties) and automorphic forms (derivatives of Eisenstein series and $L$-functions). Kudla's lecture \cite{Kudla2004} at MSRI 2001 outlined this vast program and included many inspiring discussions about history, examples and known results to date. He also proposed several central conjectures in the program:
\begin{enumerate}
\item The modularity of the generating function of special cycles in Chow groups and arithmetic Chow groups \cite[Problem 1, Problem 4]{Kudla2004};
\item The arithmetic Siegel--Weil formula, relating arithmetic volumes of the generating function to derivatives of Eisenstein series \cite[Problem 6]{Kudla2004};
\item The arithmetic inner product formula, relating heights of arithmetic theta lifts to derivatives of $L$-functions  \cite[(8.3)]{Kudla2004}.
\end{enumerate}
These conjectures were largely settled in the case of quaternionic Shimura curves over $\mathbb{Q}$, cumulating in the monograph of Kudla--Rapoport--Yang \cite{Kudla2006}. These conjectures in general for Shimura varieties of higher dimension were seemingly far from reach at that time, yet recent years have witnessed exciting advances on all of them. The goal of this article is to explain some of these recent theorems together with their applications, and mention some of the problems which remain open, providing an update to \cite{Kudla2004}.

In accordance with the theme of the summer school, we will begin with a recap of classical theory of the theta correspondence that is relevant to us (see Gan's article \cite{Gan2023} in these proceedings for more details) and develop it into the \emph{trilogy} of classical/geometric/arithmetic theta correspondences (see Table~\ref{tab:trilogy}). For certain reductive dual pairs $(G,H)$, these provide correspondences from automorphic forms for $G$ to automorphic forms/cohomology classes/algebraic cycles associated to $H$ respectively, and the last of which is the main focus of the Kudla program.  The theory of classical theta correspondences is an indispensable tool in automorphic forms and has many applications in the Langlands program. Likewise the theory of geometric/arithmetic theta correspondences is an indispensable tool in the study of cohomology/algebraic cycles of Shimura varieties and has important applications e.g., to the Hodge conjecture and Tate conjecture/the Birch--Swinnerton-Dyer conjecture and Beilinson--Bloch conjecture.

\begin{table}[h]
  \centering\setlength\extrarowheight{2pt}
  \begin{tabular}{c|c|c}
      Theta Correspondence & Lift Automorphic Forms on $G$ to & {Applications} \\\hline 
      Classical  & Automorphic Forms on $H$ & {Langlands functionality}\\\hline
      {Geometric} & {Cohomology classes on $\Sh(H)$} & 
      {Hodge conjecture} \\
      {(Kudla--Millson '80s)}& & {Tate conjecture} \\\hline
      {Arithmetic} & {Algebraic cycles on $\Sh(H)$} &
      {BSD conjecture} \\
      {(Kudla's program '90s)} & &{Beilinson-Bloch conjecture}

    \end{tabular}
  \caption{Trilogy of theta correspondences}
  \label{tab:trilogy}
\end{table}

We will mainly discuss the case when $H$ is a unitary group, for which several technical aspects are simpler than the orthogonal case  and more complete results are available, and refer to the recent exposition \cite{Li23} for some discussion of the orthogonal case including many classical examples. Needless to say, the Kudla program has grown into a vast area and many topics are not covered due to the limit on length and scope. Our discussion of the classical and geometric theta correspondence given here is only intended to highlight certain aspects directly relevant to the arithmetic version. We emphasize that a very large number of people have made fundamental contributions to the theory and we refer to Gan's article \cite{Gan2023} in these proceedings for information, references, and history. There have also been exciting recent advances on the analogue of the Kudla program over function fields and we refer to Feng--Harris's article \cite[\S4]{FH2023} in these proceedings for more details.

\section{Classical theta correspondence (recap)}\label{sec:class-theta-lift}

\subsection{Weil representation} Let $F/F_0$ be a quadratic extension of number fields. Denote by $\mathbb{A}=\mathbb{A}_{F_0}$ be the ring of adeles of $F_0$. Let $W$ be the standard split skew-hermitian space of dimension $2n$ over $F$, i.e., $W=F^{2n}$ equipped with the skew-hermitian form with matrix $w_n=\left(\begin{smallmatrix}0 & 1_n \\ -1_n & 0\end{smallmatrix}\right)$. Denote by $\mathbb{W}=W \otimes_F \mathbb{A}_F$, an skew-hermitian space of rank $2n$ over $\mathbb{A}_F$. Let $\Gn=\U(W)$, a quasi-split unitary group. Let $P=MN\subseteq G$ be the standard Siegel parabolic subgroup stabilizing the maximal isotropic subspace $F^n \oplus 0^n\subseteq W=F^{2n}$. Explicitly, under the standard basis we have  

\begin{align*}
  \Mn
  &=\left\{m(a)=\begin{pmatrix}a & 0\\0 &{}^t\bar a^{-1}\end{pmatrix}: a\in 
    \Res_{F/F_0}\GL_n
  \right\},\\
  \Nn
  &= \left\{n(b)=\begin{pmatrix} 1_n & b \\0 & 1_n\end{pmatrix}: b\in \Herm_n
  \right\}.
\end{align*} Here $a\mapsto \bar a$ is the Galois conjugation, and $\Herm_n$ is the $F_0$-subscheme  of $\Res_{F/F_0}\mathrm{Mat}_n$ given by $n\times n$ hermitian matrices $b$, i.e., ${}^t\bar b=b$.

Let $\mathbb{V}$ be an hermitian space of rank $m$ over $\mathbb{A}_F$ with hermitian form $(\ , \ )$ ($\mathbb{V}$ is not required to come from the base change of a hermitian space over $F$). Denote by $H(\mathbb{A})=\U(\mathbb{V})$, an adelic group. Let $\eta: \mathbb{A}^\times/F_0^\times\rightarrow \mathbb{C}^\times$ be the quadratic character associated to $F/F_0$. Fix $\chi: \mathbb{A}_{F}^\times\rightarrow \mathbb{C}^\times$ a character such that $\chi|_{\mathbb{A}^\times}=\eta^m$.  Such a $\chi$ determines a splitting homomorphism $$\Gn(\mathbb{A})\times H(\mathbb{A})\rightarrow \Mp(\mathbb{V} \otimes_{\mathbb{A}_F} \mathbb{W})$$ lifting the natural homomorphism $\Gn(\mathbb{A})\times H(\mathbb{A})\rightarrow \Sp(\mathbb{V} \otimes_{\mathbb{A}_F} \mathbb{W})$ (see \cite{Kudla1994}). The metaplectic group $\Mp(\mathbb{V} \otimes_{\mathbb{A}_F} \mathbb{W})$ has a distinguished representation $\omega_{\psi}=\otimes\omega_{\psi_v}$ depending on an additive character $\psi: \mathbb{A}/F_0\rightarrow \mathbb{C}^\times$. Thus a fixed choice of the pair $(\chi,\psi)$ gives rise to a \emph{Weil representation} $\omega=\omega_{\chi,\psi}=\otimes_v\omega_{\chi_v,\psi_v}$ of $\Gn(\mathbb{A})\times H(\mathbb{A})$. 
The Weil representation $\omega$ has an explicit realization on $\sS(\mathbb{V}^n)$, the space of Schwartz functions on $\mathbb{V}^n$, known as the \emph{Schr\"odinger model}: for $\varphi\in \sS(\mathbb{V}^\m)$ and $\mathbf{x}\in \mathbb{V}^n$,
\begin{align*}
\omega(m(a))\varphi(\mathbf{x})&=\chi(\det a)|\det a|_F^{\n/2}\varphi(\mathbf{x}\cdot a),&m(a)\in \Mn(\mathbb{A}),\\
\omega(n(b))\varphi(\mathbf{x})&=\psi(\tr b\,(\mathbf{x},\mathbf{x}))\varphi(\mathbf{x}),&n(b)\in \Nn(\mathbb{A}),\\
\omega(w_n)\varphi(\mathbf{x})&=\gamma_{\mathbb{V}}^n\cdot\widehat \varphi(\mathbf{x}),&w_n=\left(\begin{smallmatrix}
0  & 1_n\\
  -1_n & 0\\
\end{smallmatrix}\right),\\
\omega(h)\varphi(\mathbf{x})&=\varphi(h^{-1}\cdot\mathbf{x}),& h\in H(\mathbb{A}).
\end{align*}
Here $(\mathbf{x},\mathbf{x})=((x_i,x_j))_{1\le i,j\le n}\in\Herm_n(\mathbb{A})$ is the \emph{moment matrix} 
of $\mathbf{x}$, $\gamma_{\mathbb{V}}$ is the Weil constant (see \cite[(10.3)]{Kudla2014}), and $\widehat\varphi$ is the Fourier transform of $\varphi$ using the self-dual Haar measure on $\mathbb{V}^n$ with respect to $\psi\circ\tr_{F/F_0}$.

\subsection{Theta lifting}\label{sec:theta-lifting} Let $V$ be an $F/F_0$-hermitian space of dimension $m$
. In this case the adelic group $H(\mathbb{A})$ for $V(\mathbb{A}):=V \otimes_F \mathbb{A}_F$ agrees with the $\mathbb{A}$-points of the unitary group $H=\U(V)$ over $F_0$. We now construct theta functions using the Weil representation $\omega$ . 

Associated to $\varphi\in \sS(\VA^n)$, define the (two-variable) \emph{theta function}
\begin{equation}
  \label{eq:thetagh}
  \theta(g,h,\varphi):=\sum_{\mathbf{x}\in V^n}\omega(g,h)\varphi(\mathbf{x})=\sum_{\mathbf{x}\in V^n}\omega(g)\varphi(h^{-1}\mathbf{x}),\quad g\in \Gn(\mathbb{A}), h\in H(\mathbb{A}).
\end{equation}
 Then $\theta(g,h,\varphi)$ is as a function on $\Gn(\mathbb{A})\times H(\mathbb{A})$ is automorphic, in the sense that it is invariant under $G(F_0)\times H(F_0)$ (such invariance is a consequence of the Poisson summation formula). The construction of theta function produces the \emph{automorphic theta distribution}
 \begin{equation}
   \label{eq:autothetadis}
   \theta: \sS(\VA^n)\rightarrow\mathrm{Fun}([\Gn])) \otimes \mathrm{Fun}([H]),\quad \varphi\mapsto\theta(-,-,\varphi),
 \end{equation}
 a $G(\mathbb{A})\times H(\mathbb{A})$-equivariant distribution valued in the space of automorphic functions.

Using $\theta(g,h,\varphi)$ as an integral kernel allows one to lift automorphic forms on $G$ to automorphic forms on $H$ (and vice versa):  for an automorphic form $\phi\in \sA(\Gn(\mathbb{A}))$, define the \emph{theta lift} $\theta_\varphi(\phi)$ of $\phi$ to $H(\mathbb{A})$ by the Petersson inner product on $[G]:=G(F_0)\backslash G(\mathbb{A})$, $$\theta_\varphi(\phi)(h):=\langle \theta(-,h,\varphi),\phi\rangle_{G}=\int_{[G]}\theta(g,h,\varphi)\overline{\phi(g)}\rd g$$ if it is absolutely convergent. Here $\rd g$ is the Tamagawa measure on $G(\mathbb{A})$. If $\phi$ is cuspidal, then this integral is absolutely convergent and defines an automorphic form $\theta_\varphi(\phi)\in\sA(H(\mathbb{A}))$. Let $\pi$ be a cuspidal automorphic representation of $\Gn(\mathbb{A})$, then we obtain an $G(\mathbb{A})\times H(\mathbb{A})$-equivariant linear map $$\theta:  \sS(\VA^n) \otimes \pi^\vee\rightarrow \sA(H(\mathbb{A})),\quad (\varphi,\bar\phi)\mapsto\theta_\varphi(\phi),$$ and define the \emph{global theta lift} $\Theta_V(\pi)\subseteq \sA(H(\mathbb{A}))$ of $\pi$ to be its image, an $H(\mathbb{A})$-subrepresentation of $\sA(H(\mathbb{A}))$. The theory of global theta correspondence provides a rather complete description of $\Theta_V(\pi)$, and  we refer to Gan's article in these proceedings for more details.

\subsection{Siegel--Weil formula} Associated to $\varphi\in \sS(\VA^n)$, consider the \emph{theta integral}
\begin{equation}
  \label{eq:Igphi}
  I(g,\varphi):=\int_{[H]}\theta(g,h,\varphi)\ \rd h,
\end{equation}
  where $\rd h$ is the Tamagawa measure on $H(\mathbb{A})$
. Let $\alpha$ be the dimension of a maximal isotropic subspace of $V$.  The theta integral is absolutely convergent for all $\varphi$ if and only if the pair $(V,W)$ satisfies \emph{Weil's convergence condition}
\begin{center}
$\alpha=0$ (i.e., $V$ is anisotropic), or $\alpha>0$ and $m-\alpha>n$.
\end{center} In this case $I(g,\varphi)$ is an automorphic form on $G(\mathbb{A})$. It can be viewed as the theta lift of the identity function on $H(\mathbb{A})$, and also specializes to the weighted average of theta series within a genus class for definite hermitian forms (cf. \cite[Example 2.2.6]{Li23}). The theta integral produces a $G(\mathbb{A})$-equivariant distribution
\begin{equation}
  \label{eq:thetaintegral}
  I: \sS(\VA^n)\rightarrow\sA(G(\mathbb{A})),\quad \varphi\mapsto I(-,\varphi).
\end{equation}

There is another way of producing automorphic distributions like (\ref{eq:thetaintegral}) via Eisenstein series. For $s\in \mathbb{C}$, let $$I(s,\chi):=\Ind_{P(\mathbb{A})}^{G(\mathbb{A})}(\chi|\cdot|_F^{s+n/2})$$ be the degenerate principal series representation of $G(\mathbb{A})$, where $\Ind_{P(\mathbb{A})}^{G(\mathbb{A})}$ denotes the (unnormalized) smooth parabolic induction. Associated to $\varphi\in \sS(\VA^n)$, there is a \emph{standard Siegel--Weil section} $\Phi_\varphi(g,s)\in I(s,\chi)$ defined by $$\Phi_\varphi(g,s):=\omega(g)\varphi(0)\cdot |\det a(g)|_F^{s-s_0},$$ where  $$s_0:=\frac{m-n}{2}.$$  Here we write $g=nm(a)k$ under the Iwasawa decomposition $G(\mathbb{A})=N(\mathbb{A})M(\mathbb{A})K$ for $K$ the standard maximal open compact subgroup of $G(\mathbb{A})$, and the quantity $|\det a(g)|_F:=|\det a|_F$ is well-defined. We obtain a distribution $$\Phi(s): \sS(\VA^n)\rightarrow I(s,\chi),\quad \varphi\mapsto\Phi_\varphi(g,s),$$ and the special value $s=s_0$ is the unique value such that $\Phi(s)$ is $G(\mathbb{A})$-equivariant. 
Define the (hermitian) \Blue{Siegel Eisenstein  series} $$E(g, s,\varphi):=\sum_{\gamma\in P(F_0)\backslash G(F_0)}\Phi_\varphi(\gamma g,s),\quad g\in G(\mathbb{A}).$$  The Siegel Eisenstein series $E(g,s,\varphi)$ converges absolutely when $\Re(s)>\frac{n}{2}$. It has a meromorphic continuation\footnote{Careful readers may notice that in early works such as \cite{KR88,KR94}, the function $\varphi$ is assumed to be $K$-finite in order to obtain the meromorphic continuation. This $K$-finiteness assumption can be dropped thanks to the work of Lapid \cite{Lap2008}.} to $s\in \mathbb{C}$ and satisfies a functional equation centered at $s=0$ (see Lapid's article \cite{Lap2022} in these proceedings for the analytic theory of Eisenstein series). If $E(g,s,\varphi)$ is homomorphic at $s=s_0$, then its value at $s=s_0$ produces a $G(\mathbb{A})$-equivariant distribution
\begin{equation}
  \label{eq:eisensteinseries}
  E(s_0): \sS(\VA^n)\rightarrow\sA(G(\mathbb{A})),\quad \varphi\mapsto E(-,s_0,\varphi).
\end{equation}
The Siegel--Weil formula gives a precise identity between the two distributions (\ref{eq:thetaintegral}) and (\ref{eq:eisensteinseries}).

\begin{theorem}[Siegel--Weil formula]\label{thm:SWF} Assume that the pair $(V,W)$ satisfies Weil's convergence condition. Then $E(g,s,\varphi)$ is holomorphic at $s_0$ and $$\kappa\cdot I(g,\varphi)= E(g,s_0,\varphi),$$ where $\kappa=1/2$ if $m>n$ and $\kappa=1$ otherwise.
      \end{theorem}
This theorem was proved in Weil \cite[Theorem 5]{Weil1965} (when $m>2n$, in which case $E(g,s,\varphi)$ is also absolutely convergent at $s=s_0$), Ichino \cite[Theorem 1.1]{Ich07} (when $n<m\le 2n$) and Yamana \cite[Theorem 2.2]{Yam11} (when $m\le n$). If the Weil's convergence condition is not satisfied, one can still naturally define $I(g,\varphi)$ via regularization and it is a long effort starting with the work of Kudla--Rallis \cite{KR94} to generalize the Siegel--Weil formula  outside the convergence range and for all reductive dual pairs of classical groups. We refer to Gan--Qiu--Takeda \cite{GQT14} for the most general Siegel--Weil formula and a nice summary of the literature and history. The Siegel--Weil formula is an indispensable tool in the arithmetic theory of quadratic forms and hermitian forms (see e.g. \cite[\S1--2]{Li23} for classical examples).

\subsection{Rallis inner product formula} Piatetski-Shapiro--Rallis \cite{PR86,PR87} discovered an integral representation (the\emph{ doubling method}) of the standard $L$-function 
for cuspidal automorphic representations $\pi$ of $G(\mathbb{A})$, via integrating against a Siegel Eisenstein series on a ``doubling'' group. Combining with the doubling seesaw and the Siegel--Weil formula, one arrives at the Rallis inner product formula, which relates the Petersson inner product of theta lifts (from $G$ to $H$) and a special value of the standard $L$-function of $\pi$
. 

Consider the skew-hermitian space $W^\square=W \oplus (-W)$ of dimension $4n$ over $F$. Define $G^\square:=\U(W^\square)$, a quasi-split unitary group  of rank twice that of $\U(W)$
. Associated to the parabolic subgroup $P^\square\subseteq G^\square$ stabilizing the maximal isotropic subspace $\{(w,-w):w\in W\}\subseteq W^\square$, we have a Weil representation $\omega^\square$ of $G^\square(\mathbb{A})$ on $\sA(\VA^{2n})$. There is an isomorphism of $G(\mathbb{A})\times G(\mathbb{A})$ representations $$\delta: \omega \boxtimes (\omega^\vee \otimes \chi) \simeq\omega^\square|_{G(\mathbb{A})\times G(\mathbb{A})}$$ such that for any $\varphi_1,\varphi_2\in \sS(\VA^{n})$, we have $$\delta(\varphi_1 \otimes \overline{\varphi_2})(0)=\langle\varphi_1,\varphi_2\rangle_\omega,$$ where $\langle\ ,\ \rangle_\omega$ is the inner product on $\sS(\VA^n)$.  For any $\varphi_1,\varphi_2\in \sS(\VA^{n})$, we have a Siegel Eisenstein series $E(g,s,\delta(\varphi_1 \otimes \overline{\varphi_2}))$ on $G^\square$. 

For any $\phi_1,\phi_2\in \pi$, define the \emph{global doubling zeta integral} $$Z(s,\phi_1,\phi_2,\varphi_1,\varphi_2):=\int_{[G]\times [G]}\overline{\phi_1}(g_1)\phi_2(g_2)\cdot E((g_1,g_2),s,\delta(\varphi_1 \otimes \overline{\varphi_2}))\, \chi^{-1}(\det g_2)\, \rd g_1\rd g_2.$$ It converges absolutely when $\Re(s)\gg0$ and extends to a meromorphic function on $\mathbb{C}$.  When $\varphi_i=\otimes_v\varphi_{i,v}$ and $\phi_i=\otimes_v\phi_{i,v}$ are factorizable, the global doubling zeta integral factorizes into a product of \emph{local doubling zeta integrals} $$Z(s,\phi_1,\phi_2,\varphi_1,\varphi_2)=\prod_vZ_v(s,\phi_{1,v},\phi_{2,v}, \varphi_{1,v},\varphi_{2,v}),$$ where $$Z_v(s,\phi_{1,v},\phi_{2,v}, \varphi_{1,v},\varphi_{2,v}):=\int_{G(F_{0,v})}\overline{\langle g_v\phi_{1,v},\phi_{2,v}\rangle}_{\pi_v}\cdot \Phi_{\varphi_{1,v} \otimes \overline{\varphi_{2,v}}}((g_v,1),s)\, \rd g_v$$ converges absolutely when $\Re(s)\gg0$ and extends to a meromorphic function on $\mathbb{C}$. When all the data are normalized unramified at a finite place $v$ with $\langle\phi_{1,v},\phi_{2,v}\rangle=1$, we have $$Z_v(s,\phi_{1,v},\phi_{2,v}, \varphi_{1,v},\varphi_{2,v})=\frac{L(s+1/2,\pi_v\times\chi_v)}{b_{2n,v}(s)},$$ where $L(s+1/2, \pi_v\times \chi_v)$ is the doubling $L$-factor (see Harris--Kudla--Sweet \cite{Harris1996}, Lapid--Rallis \cite{LR05}, Yamana \cite{Yam14}) and agrees with standard (base change) $L$-factor $L(s+1/2, \mathrm{BC}(\pi_v) \otimes \chi)$ in this unramified case, and $b_{k,v}(s):=\prod_{i=1}^kL(2s+i,\eta_v^{k-i})$ is a product of Hecke $L$-factors. Define the normalized local doubling zeta integral $$Z_v^\natural(s,\phi_{1,v},\phi_{2,v}, \varphi_{1,v},\varphi_{2,v}):=\left(\frac{L(s+1/2,\pi_v\times\chi_v)}{b_{2n,v}(s)}\right)^{-1}\cdot Z_v(s,\phi_{1,v},\phi_{2,v}, \varphi_{1,v},\varphi_{2,v}),$$ then $Z_v^\natural(s,\phi_{1,v},\phi_{2,v}, \varphi_{1,v},\varphi_{2,v})=1$ for almost all $v$.

At $s=s_0$, the normalized local zeta integral evaluates to $$Z_v^\natural(s_0, \phi_{1,v},\phi_{2,v}, \varphi_{1,v},\varphi_{2,v})=\int_{G(F_{0,v})}\overline{\langle g_v\phi_{1,v},\phi_{2,v}\rangle}_{\pi_v}\cdot\langle g_v\varphi_{1,v},\varphi_{2,v}\rangle_{\omega_v} \rd g_v,$$ the integral of the product of matrix coefficients of $\pi_v$ and $\omega_v$. Thus it produces a $G(F_{0,v})\times G(F_{0,v})$-equivariant distribution $$Z_v^\natural(s_0):\sS(V_v^{2n})\rightarrow \pi_v \boxtimes (\pi_v^\vee \otimes \chi_v),\quad \varphi_{1,v}\otimes \overline{\varphi_{2,v}}\mapsto Z_v^\natural(s_0, -,-, \varphi_{1,v},\varphi_{2,v}).$$ Taking product produces a $G(\mathbb{A})\times G(\mathbb{A})$-equivariant distribution
\begin{equation}
  \label{eq:prodz}
  \prod_vZ_v^\natural(s_0): \sS(\VA^{2n})\rightarrow \pi \boxtimes (\pi^\vee \otimes \chi). 
\end{equation}

On the other hand, the Petersson inner product of  theta lifts  also defines a  $G(\mathbb{A})\times G(\mathbb{A})$-equivariant distribution
\begin{equation}
  \label{eq:thetatheta}
  \langle\theta,\theta\rangle:\sS(\VA^{2n})\rightarrow \pi \boxtimes (\pi^\vee  \otimes \chi),\quad (\varphi_1,\varphi_2)\mapsto \langle \theta_{\varphi_1}(-),\theta_{\varphi_2}(-)\rangle_{H}.
\end{equation}
 The Rallis inner product gives a precise identity between the two distributions \eqref{eq:prodz} and \eqref{eq:thetatheta}.

\begin{theorem}[Rallis inner product formula] \label{thm:rallis-inner-product}
   Assume that the pair $(V,W^\square)$ satisfies Weil's convergence condition. Let $\pi$ be a cuspidal automorphic representation of $G(\mathbb{A})$. Then for any $\phi_i=\otimes_v\phi_{i,v}\in\otimes_v\pi_{v}\simeq \pi$, $\varphi_i=\otimes_v\varphi_{i,v}\in \sS(\VA^{n})$ ($i=1,2$), $$\kappa\cdot\langle \theta_{\varphi_1}(\phi_1),\theta_{\varphi_2}(\phi_2)\rangle_{H}=\frac{L(s_0+1/2,\pi\times\chi)}{b_{2n}(s_0)}\cdot \prod_v Z_v^\natural(s_0,\phi_{1,v},\phi_{2,v}, \varphi_{1,v},\varphi_{2,v}).$$ Here $s_0=(m-2n)/2$, $\kappa=1/2$ if $m>2n$ and $\kappa=1$ otherwise, are the constants in the Siegel--Weil formula (Theorem \ref{thm:SWF}) for the pair $(V,W^\square)$.
\end{theorem}

This theorem was proved in J.-S. Li \cite{Li92}.  The formula originates from the work of Rallis \cite{Ral84} on the case of orthogonal/symplectic dual pairs
, hence its name. 
When Weil's convergence condition is not satisfied, one can still use the regularized Siegel--Weil formula to derive a regularized version of the Rallis inner product formula and again we refer to \cite{GQT14} for the most general statements.

\subsection{Theta dichotomy in the equal rank case}\label{sec:theta-dich-equal}  The Rallis inner product formula plays an important role in Rallis's program on the nonvanishing criterion for global theta lifts (\cite{Rallis1987}, cf. \cite[\S1.2]{GQT14}). We recall a special case when $m=2n$, i.e., when the two spaces $V,W$ have equal rank. In this case, $\kappa=1$ and as $m$ is even we can simply choose the splitting character $\chi$ to be the trivial character. The special point $s_0=0$ in the Siegel--Weil formula for $(V,W^\square)$  corresponds to the center of the function equation of the Eisenstein series, and the Rallis inner product formula relates the Petersson inner product of theta lifts and central $L$-values $L(1/2,\pi)$. By the Rallis inner product formula, we know that
\begin{center}
global theta lifting $\Theta_V(\pi)\ne0$ $\Longleftrightarrow$  $L(1/2,\pi)\ne0$, and $\prod_v Z_v^\natural(0)\ne0$ in \eqref{eq:prodz}.
\end{center}
The local condition $Z_v^\natural(0)\ne0$ turns out to be equivalent to that the local theta lift $\Theta_{V_v}(\pi_v)\ne0$ (cf. \cite[Proposition 3.1]{Harris1996}). Moreover, at any nonsplit place $v$ we have the \emph{theta dichotomy}: there exists a unique (up to isomorphism) local hermitian space $V_v=V_v(\pi_v)$ of rank $n$ over $F_v$ such that $Z_v^\natural(0)\ne0$ (cf. \cite[Proposition 2.6]{Liu2011}). At a finite nonsplit place $v$, we further have the \emph{epsilon dichotomy} of Harris--Kudla--Sweet \cite[Theorem 6.1]{Harris1996}, as completed by Gan--Ichino \cite[Theorem 11.1]{GI14}), pinning down exactly one of the two local hermitian spaces over $F_v$:
\begin{center}
  $Z_v^\natural(0)\ne0$ $\Longleftrightarrow$ $\varepsilon(V_v)=\omega_{\pi_v}(-1)\cdot \varepsilon(1/2,\pi_v,\psi_v)$,
\end{center}  where $\varepsilon(V_v)=\eta_v((-1)^{m(m-1)/2}\det(V_v))\in \{\pm1\}$ is the local Hasse invariant, $\varepsilon(1/2,\pi_v,\psi_v)\in\{\pm1\}$ is the central value of the doubling epsilon factor and $\omega_{\pi_v}$ is the central character of $\pi_v$. We remark that the dichotomy phenomenon involving local root numbers appeared in many earlier works, such as the classic works of Tunnell \cite{Tun1983} and Waldspurger \cite{Wal1985}, as well as the works of Prasad \cite{Pra1990} and Harris--Kudla \cite{HK1991} in the case of the triple product. We refer to Harris--Kudla--Sweet \cite{Harris1996} for the foundational theory of the theta  dichotomy for unitary groups and related history and references. 

For any cuspidal automorphic representation $\pi$ of $G(\mathbb{A})$, the theta dichotomy associates to it a \emph{unique} collection of local hermitian spaces $\{V_v=V_v(\pi_v)\}_v$ such that $\Theta_{V_v}(\pi_v)\ne0$ for all places $v$, or equivalently, a unique hermitian space $\mathbb{V}=\mathbb{V}_\pi$ of rank $n$ over $\mathbb{A}_F$ such that $\Theta_{\mathbb{V}_v}(\pi_v)\ne0$ for all places $v$, where $\mathbb{V}_v:=\mathbb{V} \otimes_{\mathbb{A}}F_{0,v}$.  Say that $\mathbb{V}$ is \emph{coherent} if $\mathbb{V}\simeq V \otimes_F \mathbb{A}_F$ for some hermitian space $V$ over $F$, and \emph{incoherent} otherwise.

Define $\varepsilon(\mathbb{V}):=\prod_{v} \varepsilon(\mathbb{V}_v)\in\{\pm1\}$. Then the Hasse principle implies that $\mathbb{V}$ is coherent if and only if $\varepsilon(\mathbb{V})=+1$.  The epsilon dichotomy implies the equality of signs $$ \varepsilon(\mathbb{V}_\pi)=\varepsilon(1/2,\pi).$$ We have two cases:
\begin{itemize}
\item If $\varepsilon(1/2,\pi)=+1$, then $\mathbb{V}_\pi$ is coherent. If $\mathbb{V}_\pi\simeq V \otimes_F \mathbb{A}_F$, then
  \begin{center}
    the global theta lift $\Theta_{V}(\pi)\ne 0$ $\Longleftrightarrow$    $L(1/2,\pi)\ne0$.
  \end{center}
Moreover $\Theta_{V'}(\pi)=0$ for all hermitian spaces $V'$ of rank $n$ over $F$ different from $V$ for local reasons.
\item If $\varepsilon(1/2,\pi)=-1$, then $\mathbb{V}_\pi$ is incoherent. The global theta lift $\Theta_{V}(\pi)=0$ for all hermitian spaces $V$ of rank $n$ over $F$ for local reasons.
\end{itemize}

In the second case there is no global theta lifting associated to the \emph{incoherent} space $\mathbb{V}=\mathbb{V}_\pi$ and $L(1/2,\pi)=0$ always, due to the sign of the functional equation.  It is natural and interesting to study the central \emph{derivative} $L'(1/2,\pi)$. The Birch and Swinnerton-Dyer conjecture and its generalization by Beilinson and Bloch suggests that the condition $L'(1/2,\pi)\ne0$ should be related to the non-triviality of \emph{algebraic cycles}. 
  When  the incoherent  space $\mathbb{V}$ is totally definite, 
 next  we will canonically associate to it a unitary Shimura variety $X$ over $F$ and use the generating function of its special cycles to define an \emph{arithmetic theta lift} $\Theta_{\mathbb{V}}(\pi)\subseteq \Ch^{n}(X)$. Here $\Ch^\mm(X)$ is the Chow group of algebraic cycles of codimension $\mm$ on $X$ modulo rational equivalence. One of the goals of the Kudla program on arithmetic theta lifting is to establish an analogous criterion (cf. Theorem \ref{thm:AIPF})
 \begin{equation}
   \label{eq:goal}
\text{       the arithmetic theta lift $\Theta_{\mathbb{V}}(\pi)\ne 0$ $\stackrel{?}{\Longleftrightarrow}$   $L'(1/2,\pi)\ne0$.}
\end{equation}

\section{Geometric theta correspondence}

\subsection{Unitary Shimura varieties} \label{sec:unit-shim-vari}From now on we assume that $F/F_0$ is a CM extension of a totally real number field. 

As in \S\ref{sec:class-theta-lift}, $V$ denotes a hermitian space over $F$ of rank $\n$ and  $H=\U(V)$. We fix an embedding $\sigma: F\hookrightarrow \mathbb{C}$ and view $F$ (resp. $F_0$)
as a subfield of $\mathbb{C}$ (resp. $\mathbb{R}$)
. Say $V$ is \emph{standard indefinite} if $V$ has signature $(\n-1,1)$ at the  real place of $F_0$ induced by $\sigma$
, and signature $(\n,0)$ at all other real places.  When $V$ is standard indefinite, there is a system of unitary Shimura varieties $X=\{X_K\}$ indexed by neat open compact subgroup $K\subseteq H(\mathbb{A}_f)$. Each $X_K$ is a smooth quasi-projective scheme of dimension $\n-1$ over $F\subseteq \mathbb{C}$, and is projective when $V$ is anisotropic (e.g., when $F_0\ne \mathbb{Q}$, by the signature condition). It has complex uniformization $$X_K(\mathbb{C})=H(F_0)\backslash[ \mathbb{D}\times H(\mathbb{A}_f)/K],$$  where $\mathbb{D}$ is the hermitian symmetric domain associated to $\U(V_\infty)$ given by the space of negative complex lines in $V \otimes_F\mathbb{C}$. We have isomorphisms $$\mathbb{D}\simeq\{z\in \mathbb{C}^{\n-1}: |z|< 1\}\simeq \frac{\UU(\n-1,1)}{\UU(\n-1)\times\UU(1)}.$$  In particular, $X_K$ can be written as a union of arithmetic quotients of complex balls.  

As in \S\ref{sec:class-theta-lift},  $\mathbb{V}$ denotes an  incoherent hermitian space of rank $\n$ over $\mathbb{A}_F$. Say $\mathbb{V}$ is \emph{totally definite} if $\mathbb{V}$ has signature $(\n, 0)$ at all real places. If $\mathbb{V}$ is totally definite, then for any embedding $\sigma: F\hookrightarrow \mathbb{C}$, we have a unique standard indefinite hermitian space $V$, depending on $\sigma$, such that $V_v$ has signature  $(\n-1,1)$ at the  real place of $F_0$ induced by $\sigma$, and $\mathbb{V}_v\simeq V_v$ at all other places of $F_0$. By the theory of conjugation of Shimura varieties, the Shimura variety $X_K$ associated to varying $V$ for varying choices of $\sigma$ are all conjugate, and thus can be intrinsically defined over $F$ (without being viewed  as a subfield of $\mathbb{C}$). In other words, for any totally definite incoherent hermitian space $\mathbb{V}$ over $\mathbb{A}_F$, we obtain a system of unitary Shimura varieties $X=\{X_K\}$ canonically defined over $F$ (cf. \cite{Zha19,Gro21}).

From the above discussion the following dichotomy picture emerges:
\begin{itemize}
\item when studying the \emph{geometric} invariants of $X_K$ (over the algebraically closed field $\mathbb{C}$), a choice of the embedding $\sigma:F\hookrightarrow \mathbb{C}$ is involved. The \emph{coherent} space $V(\mathbb{A})$ associated to $V$ should play a canonical role and special \emph{values} of analytic quantities ought to appear.
\item when studying the \emph{arithmetic} invariants of $X_K$ (over the number field $F$), no choice of the embedding $\sigma:F\hookrightarrow \mathbb{C}$ is involved and the \emph{incoherent} space $\mathbb{V}$ should play a canonical role  and special \emph{derivatives} of analytic quantities ought to appear.
\end{itemize}

\begin{remark}
  By the Langlands philosophy, the motivic $L$-function associated to the \etale cohomology of $X_K$ should be factorized into a product of automorphic $L$-functions for automorphic representations $\pi$ of $H(\mathbb{A})$ (see Morel's article \cite{Morel2023} in these proceedings for more details). When $V$ is standard indefinite, the $L$-function appearing should be the \emph{standard} $L$-function of $\pi$. This suggests the terminology and its relevance for our goal \eqref{eq:goal}. When $V$ has more general signature combinations, for the corresponding Shimura varieties one expects to see Langlands $L$-functions  associated to more complicated 
representations of the dual group of $H=\U(V)$ (cf. \cite[Example 3.7 (2)]{Morel2023}).  
\end{remark}

\begin{remark}
  We remark that $X_K$ is a Shimura variety \Blue{of abelian type} (rather than of PEL or Hodge type). Unlike Shimura varieties of  PEL type associated to unitary similitude groups, it lacks a good moduli description in terms of abelian varieties with additional structures and thus it is technically more difficult to study. Nevertheless, its \'etale cohomology and $L$-function will be computed in terms of automorphic forms in the forthcoming work of Kisin--Shin--Zhu \cite{Kisina}, with the help of the endoscopic classification for unitary groups due to Mok \cite{Mok15} and Kaletha--Minguez--Shin--White \cite{KMSW14}. 
\end{remark}

\subsection{Special cycles} Assume that $V$ is standard indefinite and let $\mathbb{V}$ be the associated totally definite incoherent space. Let $\VAf:=V \otimes_{F_0}\mathbb{A}_f\simeq \mathbb{V}\otimes_{\mathbb{A}}\mathbb{A}_f$. For any $y\in V$ with $(y,y)>0$, i.e., with totally positive norm, its orthogonal complement $V_y\subseteq V$ is a standard indefinite hermitian space  rank $\n-1$ over $F$. Let $H_y:=\U(V_y)$, a subgroup of $H=\U(V)$ and $X_y$ be the system of unitary Shimura varieties associated to $H_y$. We define the \emph{special divisor} $Z(y)_K$ to be the Shimura subvariety $$Z(y)_K:=(X_y)_{K\cap H_y(\mathbb{A}_f)}\rightarrow X_K.$$ More generally, for any $x\in \VAf$ with $(x,x)\in F_{>0}$, there exists $y\in V$ and $h\in H(\mathbb{A}_f)$ such that $y=hx$. Define the \Blue{special divisor} $Z(x)_K$ to be the Hecke translate of a Shimura subvariety $$Z(x)_K:=(X_y)_{hKh^{-1}\cap H_y(\mathbb{A}_f)}\rightarrow X_{hKh^{-1}}\xrightarrow{\cdot h} X_K.$$

For any $\mm\le\dim X_K$ and any $\mathbf{x}=(x_1,\ldots, x_\mm)\in \VAf^\mm$ with $(\mathbf{x},\mathbf{x})\in \Herm_{\mm}(F_0)_{>0}$, define the \Blue{special cycle} (of codimension $\mm$)
\begin{equation}
  \label{eq:Zbx}
  Z(\mathbf{x})_K=Z(x_1)_K\cap \cdots \cap Z(x_\mm)_K\rightarrow X_K,
\end{equation}
  here $\cap$ denotes the fiber product over $X_K$, whose image cycle defines an algebraic cycle of codimension $\mm$ on $X_K$. It only depends on the $F$-span $V_\mathbf{x}$ of $\{x_1,\ldots,x_\mm\}$ in $\VAf^\mm$ and we write $Z(V_\mathbf{x})_K:=Z(\mathbf{x})_K$.

More generally, when $(\mathbf{x},\mathbf{x})\in \Herm_{\mm}(F_0)_{\ge0}$ but is singular, the intersection (\ref{eq:Zbx}) is improper, i.e., has the wrong codimension. Let $\mathcal{L}_K$ be the tautological line bundle on $X_K$, with complex uniformization $$\mathcal{L}_K(\mathbb{C})=H(F_0)\backslash [\mathcal{L}\times H(\mathbb{A}_f)/K],$$ where $\mathcal{L}$ is the tautological line bundle on $\mathbb{D}\subseteq\mathbb{P}(V \otimes_F \mathbb{C})$. The tautological line bundle naturally appears when computing improper intersections: for example if $(x,x)>0$, then the excess intersection formula implies that $$Z(x)_K\cdot Z(x)_K=Z(x)_K\cdot c_1(\mathcal{L}_K^\vee)\in \Ch^2(X_K).$$  Here $c_1(\mathcal{L}_K^\vee)\in \Ch^1(X_K)$ is the first Chern class of the dual line bundle of $\mathcal{L}_K$. This motivates us to define  $$Z(\mathbf{x})_K:=Z(V_\mathbf{x})_K\cdot c_1(\mathcal{L}_K^\vee)^{\mm-\dim_FV_\mathbf{x}}\in \Ch^\mm(X_K),$$ which is an element in the Chow group of correct codimension. In particular, when $\mathbf{x}=0\in \VAf^\mm$, we have $Z(\mathbf{x})_K= c_1(\mathcal{L}_K^\vee)^{\mm}\in \Ch^\mm(X_K)$.

For a $K$-invariant Schwartz function $\varphi\in \sS(\VAf^\mm)^K$ and $T\in \Herm_{\mm}(F_0)_{\ge 0}$,  define the \Blue{weighted special cycle} $$Z(T, \varphi)_K=\sum_{\mathbf{x}\in K\backslash \VAf^\mm\atop (\mathbf{x},\mathbf{x})=T}\varphi(\mathbf{x}) Z(\mathbf{x})_K\in \CHXK.$$

\subsection{Kudla's generating function (arithmetic theta function)}\label{sec:kudl-gener-funct}

Define \Blue{Kudla's generating function of special cycles} (of codimension $\mm$)
\begin{equation}
  \label{eq:kudlagen}
  Z(\sz,\varphi)_K=\sum_{T\in\Herm_\mm(F_0)_{\ge0}}Z(T,\varphi)_K\cdot q^T,
\end{equation}
as a formal generating function valued in $\CHXK$, where $$\sz\in \mathcal{H}_\mm=\{\sx+i\sy: \sx\in\Herm_\mm(F_{0,\infty}),\ \sy\in\Herm_\mm(F_{0,\infty})_{>0}\}$$ lies in the hermitian half space and $q^T:=\prod_{v\mid \infty}e^{2\pi i \tr T\sz_v}.$ It formally resembles the Fourier expansion of a holomorphic hermitian modular form on $\mathcal{H}_n$. In fact its modularity is the content of Kudla's modularity conjecture (see Conjecture \ref{conj:modularity}).

More adelically, define $$Z(g,\varphi)_K:=\sum_{T\in\Herm_\mm(F_0)_{\ge0}} Z(T,\omega_f(g_f)\varphi)_K\cdot \omega_\infty(g_\infty)\varphi_\infty(T),\quad g\in G(\mathbb{A})$$ as a formal sum valued in $\CHXK$. Here $\varphi_\infty\in\sS(\mathbb{V}_\infty^\mm)$ is the standard Gaussian function $\varphi_\infty(\mathbf{x}):=\prod_v e^{-2\pi \tr(\mathbf{x},\mathbf{x})}$ and $\omega_\infty(g_\infty)\varphi_\infty(T)$ makes sense as $\omega_\infty(g_\infty)\varphi_\infty$ factors through the moment map $\mathbf{x}\mapsto(\mathbf{x},\mathbf{x})$. It is the adelization of \eqref{eq:kudlagen} and agrees with the formal Fourier expansion of
\begin{equation}
  \label{eq:adelicgen}
  Z(g,\varphi)_K=\sum_{\mathbf{x}\in K\backslash \VAf^\mm}\omega(g)(\varphi \otimes \varphi_\infty)(\mathbf{x})\cdot Z(\mathbf{x})_K,
\end{equation} where for $\mathbf{x}\in \mathbb{V}_f^\mm$ we interpret $\varphi_\infty(\mathbf{x})$ as $\varphi_\infty((\mathbf{x},\mathbf{x}))$ if $(\mathbf{x},\mathbf{x})\in\Herm_\mm(F_0)_{\ge0}$ and 0 otherwise. Moreover $Z(g,\varphi)_K$ is compatible under pullback by natural projection morphisms when varying $K\subseteq H(\mathbb{A}_f)$ and thus defines a formal sum $$Z(g,\varphi):=(Z(g,\varphi)_K)_K$$ valued in $\CHX:=\varinjlim_{K\subseteq H(\mathbb{A}_f)}\CHXK$.

Notice the analogy between the classical theta function \eqref{eq:thetagh} and Kudla's generating function \eqref{eq:adelicgen}, except two crucial modifications:
\begin{enumerate}
\item the automorphic forms space $\sA(H(\mathbb{A}))$ (second variable) is replaced by $\Ch^\mm(X)$ for the system of Shimura varieties $X$ associated to $H$. 
\item the holomorphy of $Z(g,\varphi)$ forces us to \emph{fix} $\varphi_\infty$ to be the Gaussian function, and $\varphi_\infty$ lives on the totally definite incoherent space $\mathbb{V}$ rather than the standard indefinite space $V$ (which matches the dichotomy philosophy as discussed in \S\ref{sec:unit-shim-vari}).
\end{enumerate} In this way one should view Kudla's generating function  as an \emph{arithmetic theta function}.

\subsection{Geometric modularity} \label{sec:geometric-modularity} We can extract geometric invariants of an element $Z\in \Ch^\mm(X_K)$ by taking its Betti cohomology class $[Z]\in \H^{2\mm}(X_K(\mathbb{C}), \mathbb{Z})$ of the complex manifold $X_K(\mathbb{C})$. In particular, we obtain from the arithmetic theta function $Z(g,\varphi)_K$ a  \emph{geometric theta function} $[Z(g,\varphi)_K]$ valued in $\HXK$. Its Fourier coefficients encodes the information about the geometric intersection numbers of special cycles. The classical theorem of Kudla--Millson shows that this geometric theta function is indeed \emph{modular}. In other words, there are many \emph{hidden symmetry and relations} between these geometric invariants of special cycles.

  More precisely, denote by $\Ah\subset\sA(G(\mathbb{A}))$ the adelization of the space of  holomorphic hermitian modular forms on $\mathcal{H}_n$ of parallel weight $m$, that is, the space spanned by automorphic forms $\phi$ on $G(\mathbb{A})$ such that $\phi_\infty$ is in the minimal $K$-type of a discrete representation of weight $ (( m-\mathfrak{k}^\chi ) /2,( m+\mathfrak{k}^\chi ) /2)$ of $G ( F_{0,\infty} )$   as defined on \cite[Page 854]{Liu2011}.

\begin{theorem}[Geometric modularity]
   The formal generating function $[Z(g,\varphi)_K]$ converges absolutely and defines an element in $\Ah \otimes \HXK$.
\end{theorem}

This theorem is proved in Kudla--Millson \cite{KM90}. In fact \cite{KM90} proves a much more general theorem, applicable to the generating function of special cohomology classes for locally symmetric spaces associated to any $\U(p,q)$ or $\O(p,q)$. The proof replies on the \emph{Kudla--Millson Schwartz forms} (\cite{Kudla1986,KM87}) $$\varphi_{\KM,v_0}\in \sS(V_{v_0}^\mm) \otimes \Omega^{\mm,\mm}(\mathbb{D}),$$ where $v_0$ is the real place of $F_0$ induced by the fixed embedding $\sigma:F\hookrightarrow \mathbb{C}$
, and $\Omega^{a,b}(\mathbb{D})$ is the space of smooth differential forms on $\mathbb{D}$ of type $(a,b)$. The Schwartz form $\varphi_{\KM,v_0}$ is $H_{v_0}(\mathbb{R})$-invariant and \emph{closed} at any $\mathbf{x}\in V_{v_0}^\mm$.  Define $$\widetilde\varphi_\infty=\varphi_{\KM,v_0} \otimes \bigotimes_{v\mid \infty, v\ne v_0}\varphi_{v}\in  \sS(V_\infty^\mm)\otimes \Omega^{\mm,\mm}(\mathbb{D}),$$ where $\varphi_{v}\in \sS(V_{v}^\mm)$ is the Gaussian function. Define
\begin{equation}
  \label{eq:Schwartzform}
  \varphit:=\varphi \otimes \wit\varphi_\infty\in \sS(V(\mathbb{A})^\mm)\otimes \Omega^{\mm,\mm}(\mathbb{D})
\end{equation}
 and the \emph{Kudla--Millson theta function} $$\theta_\KM(g,h,\varphi):=\sum_{\mathbf{x}\in V^n}\omega(g)\varphit(h^{-1}\mathbf{x}),\quad g\in \Gn(\mathbb{A}), h\in H(\mathbb{A}_f),$$ which gives a closed $(\mm,\mm)$-form on $X_K(\mathbb{C})$ at any $g\in G(\mathbb{A})$. By the Poisson summation formula one can prove that  $\theta_\KM(g,h,\varphi)$ defines a (nonholomorphic) automorphic form valued in closed  $(\mm,\mm)$-forms on $X_K(\mathbb{C})$. \cite{KM90} further proves that it represents the (holomorphic) geometric theta series $[Z(g,\varphi)_K]$ in $\HXK$ (in particular, the nonholomorphic terms in $\theta_\KM(g,h,\varphi)$ are exact forms) and obtains the theorem.

\begin{remark}
  The remarkable discovery that generating function involving intersection numbers of algebraic cycles are modular originates from the work of Hirzebruch--Zagier \cite{Hirzebruch1976} on Hilbert modular surfaces, and is one of the inspirations for Kudla's work (cf. the introduction of \cite{KM90,Kudla1997}). 
    
\end{remark}

\subsection{Geometric theta lifting}\label{sec:geom-theta-lift}
Analogous to \S\ref{sec:theta-lifting}, using $[Z(g,\varphi)]$ as an integral kernel allows one to lift automorphic forms on $G$ to cohomology classes on $X(\mathbb{C})$. For $\phi\in \Ah$, define the \emph{geometric theta lift} or \emph{Kudla--Millson lift} $\theta^\KM_\varphi(\phi)$  to be the Petersson inner product $$\theta^\KM_\varphi(\phi)_K:=\langle [Z(g,\varphi)_K],\phi\rangle_G=\int_{[G]}[Z(g,\varphi)]\overline{\phi(g)}\rd g\in \HXK.$$  When varying $K\subseteq H(\mathbb{A}_f)$ it defines a class $$\theta_\varphi^\KM(\phi):=(\theta_\varphi^\KM(\phi)_K)_K\in \HX:=\varinjlim_{K\subseteq H(\mathbb{A}_f)}\HXK.$$ Let $\pi$ be a cuspidal automorphic representation of $\Gn(\mathbb{A})$. Assume that $\pi\cap\Ah\ne0$, which forces that $\pi_\infty$ is a holomorphic discrete series of a particular weight (cf. \cite[p.852]{Liu2011}). Then we obtain an $G(\mathbb{A}_f)\times H(\mathbb{A}_f)$-equivariant linear map $$\theta^\KM:  \sS(\VAf^n) \otimes \pi^\vee\rightarrow\HX,\quad (\varphi,\bar\phi)\mapsto\theta^\KM_\varphi(\phi).$$ Define the \emph{geometric theta lift} $\Theta^\KM_V(\pi)\subseteq \HX$ of $\pi$ to be its image (here $\theta^\KM_\varphi(\phi)$ is understood to be 0 if $\phi\not\in\Ah$).

\subsection{Geometric Siegel--Weil formula} To relate the geometric theta series to Eisenstein series, we need to extract numerical invariants from cohomology classes. To that end, assume that $V$ is anisotropic, thus $X_K$ is projective and we have a degree map $\deg: \H^{2\dim X_K}(X_K,\mathbb{C})\rightarrow \mathbb{C}$. 
For any $\mm\le\dim X_K=m-1$, define the \emph{geometric volume} $$\vol: \HXK\rightarrow \mathbb{C}, \quad[Z]\mapsto\deg([Z]\cup[c_1(\mathcal{L}_K^\vee)]^{\dim X_K-\mm }).$$ In particular when $n=0$ we obtain the geometric volume  $\vol([X_K])$ of the Shimura variety $X_K$. Define the \emph{normalized geometric volume} $$\voln: \HXK\rightarrow \mathbb{C}, \quad[Z]\mapsto \frac{\vol([Z])}{\vol([X_K])/2}.$$  The Haar measure on $H(\mathbb{A}_f)$ such that $K$ has volume $(\vol([X_K])/2)^{-1}$  can be viewed as an analogue of the Tamagawa measure on $H(\mathbb{A})$ (cf. \cite[Footnote 11]{LL2020}), hence the normalization. Then $\voln([Z(g,\varphi)_K])$ is independent of the choice of $K$ and can be viewed as a geometric analogue of the theta integral \eqref{eq:Igphi} and  produces a $G(\mathbb{A}_f)$-equivariant distribution analogous to \eqref{eq:thetaintegral}
\begin{equation}
  \label{eq:voln}
  \voln: \sS(\VAf^n)\rightarrow\sA(G(\mathbb{A})),\quad \varphi\mapsto \voln[Z(-,\varphi)].
\end{equation}

On the other hand,  for any $\varphi\in \sS(\mathbb{V}_f^\m)^K$, the Schwartz form $\varphit$ in \eqref{eq:Schwartzform} 
gives an element $\varphic\in\sS(V(\mathbb{A})^\mm)) \otimes \Omega^{(0,0)}(\mathbb{D})$ such that $$\varphit\wedge \Omega^{\dim X-\mm}=\varphic\cdot\Omega^{\dim X},$$ where $\Omega\in\Omega^{(1,1)}(\mathbb{D})$ is the first Cherm form of $\mathcal{L}^\vee$. Evaluation of $\varphic$ at the base point of $\mathbb{D}$ gives a Schwartz function in $\sS(V(\mathbb{A})^\mm)$, which we still denote by $\varphic$ by abuse of notation. Hence we obtain a coherent Siegel Eisenstein series $E(g,s,\varphic)$ on $G(\mathbb{A})$.

\begin{theorem}[Geometric Siegel--Weil formula]\label{thm:geoSWF} Assume that $V$ is anisotropic. Assume that $\mm\le\dim X_K=\n-1$. Then for any $\varphi\in \sS(\mathbb{V}_f^\m)$, the following identity holds
  $$\kappa\cdot \voln([Z(g,\varphi)])= E(g,s_0,\varphic).$$ Here $s_0=(m-n)/2$, $\kappa=1/2$ are the constants in the Siegel--Weil formula (Theorem \ref{thm:SWF}) for the pair $(V,W)$. 
\end{theorem}

This is \cite[(4.4)]{Kudla2004} (see also \cite[Theorem 4.23]{Kud03}) for orthogonal Shimura varieties. We refer to \cite[\S2.2]{Dun22} for an exposition of the proof for the unitary Shimura variety $X_K$. As a consequence, the geometric volumes of the special cycles $Z(T, \varphi)_K$ are related to the Fourier coefficients of the Eisenstein series $E(g,s_0,\varphic).$. The geometric Siegel--Weil formula holds more generally for non-projective $X_K$ under Weil's convergence condition, although the geometric volume lacks a cohomological interpretation as above (see \cite[Theorem 4.1]{Kudla2004}).

\subsection{Geometric inner product formula} To finish the geometric story, we introduce a geometric analogue of the Petersson inner product on $[H]$. Further assume that $2\mm\le\dim X_K$. Define the (normalized) \emph{geometric inner product} $$\langle \ ,\ \rangle_{X_K(\mathbb{C})}: \HXK\times \HXK\rightarrow \mathbb{C},\quad ([Z_1],[Z_2])\mapsto \voln([Z_1]\cup\overline{[Z_2]}).$$ It is again compatible when varying $K$ and thus gives an inner product $\langle\ , \ \rangle_{X(\mathbb{C})}$ on $\HX$. Combining the geometric Siegel--Weil formula and the Rallis inner product formula we obtain the following.

\begin{theorem}[Geometric inner product formula] \label{thm:geom-inner-prod}
   Assume that $V$ is anisotropic. Assume that $2\mm\le\dim X_K=\n-1$. Let $\pi$ be a cuspidal automorphic representation of $G(\mathbb{A})$ such that $\pi\cap\Ah\ne0$. Then for any $\phi_i=\otimes_v\phi_{i,v}\in \pi\cap\Ah$, $\varphi_i=\otimes_v\varphi_{i,v}\in \sS(\VAf^{n})$ ($i=1,2$), $$\kappa\cdot\langle\theta^\KM_{\varphi_1}(\phi_1),\theta^\KM_{\varphi_2}(\phi_2)\rangle_{X(\mathbb{C})}=\frac{L(s_0+1/2,\pi\times\chi)}{b_{2n}(s_0)}\cdot \prod_v Z_v^\natural(s_0,\phi_{1,v},\phi_{2,v}, \varphic_{1,v},\varphic_{2,v}).$$ Here $s_0=(m-2n)/2$, $\kappa=1/2$ are the constants in the Siegel--Weil formula (Theorem \ref{thm:SWF}) for the pair $(V,W^\square)$.
 \end{theorem}

 \begin{example}
 In the special case $2\mm=m-1$, each $\theta^\KM_{\varphi_i}(\phi_i)$ is the cohomology class of a middle dimensional cycle on $X(\mathbb{C})$ and the geometric inner product relates their geometric intersection number to the near central value $L(1,\pi\times\chi)$ at $s_0=1/2$.   
 \end{example}

 Kudla--Millson's theory of geometric theta correspondence \cite{KM90}, as extended by Funke--Millson to nontrivial coefficients \cite{FM06} and compactifications of non-compact $X_K$ (e.g. \cite{FM14}), have many applications to the cohomology of Shimura varieties and more general locally symmetric spaces. For example, Bergeron--Millson--Moeglin \cite{BMM16} proved the Hodge conjecture and the Tate conjecture for the arithmetic ball quotients $X_K$, in codimension $\le \frac{1}{3}\dim X_K$ or $\ge \frac{2}{3}\dim X_K$
 , and geometric theta lifting is a key ingredient in the proof to generate many Hodge/Tate classes using special cycles in these degrees far away from the middle degree. Analogous to the classical theory, 
 the geometric inner product formula and its variants (e.g. \cite[Theorem 1.1]{BF10}) are useful to prove nonvanishing results on geometric theta lifting.

\section{Arithmetic theta correspondence}\label{sec:arithm-theta-lift}

\subsection{Arithmetic modularity conjecture and Arithmetic theta lifting}\label{sec:arithm-modul-conj}

The modularity of classical and geometric theta functions motivates Kudla
s  arithmetic modularity conjecture \cite[Problem 1]{Kudla2004}.

  \begin{conjecture}[Arithmetic modularity] \label{conj:modularity}
   The formal generating function $Z(g,\varphi)_K$ converges absolutely and defines an element in $\Ah \otimes \CHXK$.
 \end{conjecture}

The formulation in the unitary case can be found in Liu \cite{Liu2011}, who also proved the case $\m=1$ and reduce the $\m>1$ case to the convergence. Recently Xia \cite{Xia22} proved the desired convergence when $F=\mathbb{Q}(\sqrt{-d})$ for $d=1,2,3,7,11$, and thus established Conjecture~\ref{conj:modularity} in these cases.

 \begin{remark}
 Kudla's arithmetic modularity conjecture was originally formulated for orthogonal Shimura varieties over $\mathbb{Q}$ (\cite{Kudla1997},\cite[Problem 1]{Kudla2004}). In this case, Borcherds \cite{Bor99} proved the conjecture for the divisor case $\m=1$. The special case of Heegner points on modular curves dates back to the classical work of Gross--Kohnen--Zagier \cite{GKZ87} (see \cite[Example 6.4.1]{Li23} for an explicit example). Zhang \cite{Zha09} proved the modularity for general $\m$ assuming the absolute convergence of the series. Bruinier--Westerholt-Raum \cite{BW15} proved the desired convergence and hence established Kudla's modularity conjecture for orthogonal Shimura varieties over $\mathbb{Q}$.     For orthogonal Shimura varieties over totally real fields, Yuan--Zhang--Zhang \cite{YZZ09} proved the modularity for $\m=1$ (see also Bruinier \cite{Bru12} for a different proof) and reduce the $\m>1$ case to the convergence. More recently, Bruinier--Zemel \cite{BZ22} has extended the modularity to toroidal compactifications of orthogonal Shimura varieties when $\m=1$ (see also Engel--Greer--Tayou \cite{EGT23} for a different proof). For more general classes of orthogonal and unitary Shimura varieties (indefinite at possibly multiple archimedean places), the modularity conjecture is proved assuming the conjectural injectivity of the Abel--Jacobi maps  by  Kudla \cite{Kud2021} and Maeda \cite{Mae2021,Mae2022}.
\end{remark}

Assume Conjecture \ref{conj:modularity}.  Analogous to \eqref{eq:autothetadis}, we obtain an \emph{arithmetic theta distribution}  $$Z: \sS(\mathbb{V}_f^\mm)\rightarrow \sA(G(\mathbb{A})) \otimes\CHX,\quad \varphi\mapsto Z(-,\varphi).$$ It is $G(\mathbb{A}_f)\times H(\mathbb{A}_f)$-equivariant, where $H(\mathbb{A}_f)$ acts on $\CHX$ via the Hecke correspondences. Analogous to \S\ref{sec:geom-theta-lift}, using $Z(g,\varphi)$ as an integral kernel allows one to lift automorphic forms on $G$ to algebraic cycles on $X$. For $\phi\in \Ah$, define the \emph{arithmetic theta lift} $\Theta_\varphi(\phi)$  to be the Petersson inner product $$\Theta_\varphi(\phi)_K:=\langle Z(g,\varphi)_K,\phi\rangle_G=\int_{[G]}Z(g,\varphi)\overline{\phi(g)}\rd g\in \CHXK.$$  When varying $K\subseteq H(\mathbb{A}_f)$ it defines a class $$\Theta_\varphi(\phi):=(\Theta_\varphi(\phi)_K)_K\in \CHX.$$ Let $\pi$ be a cuspidal automorphic representation of $\Gn(\mathbb{A})$ and assume that $\pi\cap\Ah\ne0$. Then we obtain an $G(\mathbb{A}_f)\times H(\mathbb{A}_f)$-equivariant linear map $$\Theta:  \sS(\VAf^n) \otimes \pi^\vee\rightarrow\CHX,\quad (\varphi,\bar\phi)\mapsto\Theta_\varphi(\phi),$$ and define the \emph{arithmetic theta lift} $\Theta_\mathbb{V}(\pi)\subseteq \CHX$ of $\pi$ to be its image. Again here $\Theta_\varphi(\phi)$ is understood to be 0 if $\phi\not\in\Ah$. In particular, $\Theta$ can be viewed as an element $$\Theta\in\Hom_{H(\mathbb{A}_f)}((\sS(\VAf^n)\otimes \pi_f^\vee)_{G(\mathbb{A}_f)},\CHX).$$ Notice that $(\sS(\VAf^n)\otimes \pi_f^\vee)_{G(\mathbb{A}_f)}$ is nothing but the classical theta lift $\Theta(\pi_f):=\otimes_{v\nmid\infty}\Theta_{V_v}(\pi_v)$ of $\pi_f$, thus we may view the arithmetic theta lift as an element of the $\Theta(\pi_f)$-isotypic part of $\CHX$, $$\Theta\in\Hom_{H(\mathbb{A}_f)}(\Theta(\pi_f),\CHX).$$

\subsection{Special cycles on integral models}
Kudla \cite[Problem 4]{Kudla2004} also proposed the modularity problem in the \emph{arithmetic Chow group} $\aCh^\m(\mathcal{X}_K)$, where \emph{arithmetic intersection theory} takes place  (see \cite{GS90,BGKK07} and also \cite{Sou92}), of a suitable (compactified) regular integral model $\mathcal{X}_K$ (of a variant) of $X_K$. For the purpose of this article, it suffices to know that elements in $\aCh^\m(\mathcal{X}_K)$ can be represented by $\widehat{\mathcal{Z}}=(\mathcal{Z}, (g_{\mathcal{Z},\sigma})_{\sigma:F\hookrightarrow \mathbb{C}})$, where
      \begin{enumerate}
      \item $\mathcal{Z}$ is codimension $\m$ cycle on $\XK$.
      \item $g_{\mathcal{Z},\sigma}$ is a Green current for $\mathcal{Z}_\sigma(\mathbb{C})$.
      \end{enumerate}
The data $\mathcal{Z}$ (resp. $(g_{\mathcal{Z},\sigma})_{\sigma:F\hookrightarrow \mathbb{C}}$) encodes information at finite places (resp. infinite places) for arithmetic intersection theory.

The problem \cite[Problem 4]{Kudla2004} seeks to define canonically an explicit arithmetic generating function  $\widehat{\mathcal{Z}}(\sz,\varphi)$ valued in $\aCh^\m(\mathcal{X})_\mathbb{C}$ which lifts $Z(\sz,\varphi)$ under the restriction map $$\aCh^\m(\mathcal{X})\rightarrow\Ch^\m(X),$$ and such that $\widehat{\mathcal{Z}}(\sz, \varphi)$ is modular.

To define the integral model and special cycles on it, it is more convenient to work with a related unitary Shimura variety with an explicit moduli interpretation after \cite{Kudla2014,Bruinier2017,Rapoport2017}. Define a torus $Z^\mathbb{Q}=\{z\in \Res_{F/\mathbb{Q}}\mathbb{G}_m: \mathrm{Nm}_{F/F_0}(z)\in \mathbb{G}_m\}$. Fix  a CM type $\Phi\subseteq \Hom(F, \overline{\mathbb{Q}})$ of $F$. Then associated to $\wit H\coloneqq Z^\mathbb{Q}\times \Res_{F_0/\mathbb{Q}}H$ there is a natural Shimura datum $(\wit H,\{h_{\wit H}\})$ of PEL type (\cite[\S11.1]{LZ}). Assume that $K_{Z^\mathbb{Q}}\subseteq Z^\mathbb{Q}(\mathbb{A}_f)$ is the unique maximal open compact subgroup. Then the associated Shimura variety $\Sh_{\KG}=\Sh_{K_{Z^\mathbb{Q}}\times K}(\wit H,\{h_{\wit H}\})$ is of dimension $n-1$ and has a canonical model over its reflex field $E$. Moreover. $\Sh_K$ can be identified as the product of the base change $(X_K)_E$ and a 0-dimensional Shimura variety of PEL type 
(\cite[Lemma 5.2]{LL2020}).

 Assume that $K=\prod_{v\nmid \infty} K_{v}\subseteq H(\mathbb{A}_f)$ and $K_{v}\subseteq H(F_{0,v})$ is given by
\begin{itemize}
\item the stabilizer of a self-dual or almost self-dual lattice $\Lambda_v\subseteq V_v$ if $v$ is inert in $F$,
\item the stabilizer of a self-dual lattice $\Lambda_v\subseteq V_v$ if $v$ is ramified in $F$,
\item a principal congruence subgroup of $H_v(F_{0,v})\simeq\GL_n(F_{0,v})$ if $v$ is split in $F$.
\end{itemize}

Let $\mathcal{V}_\mathrm{ram}$ (resp. $\mathcal{V}_\mathrm{asd}$) be the set of finite places $v$ of $F_0$ such that $v$ is ramified in $F$ (resp. $v$ is inert in $F$ and $\Lambda_v$ is almost self-dual). 
Further assume that all places of $E$ above $\mathcal{V}_\mathrm{ram}\cup\mathcal{V}_\mathrm{asd}$ are unramified over $F$
. Then we obtain a global regular integral model $\XK$ of $\Sh_K$ over $O_E$ after Rapoport--Smithling--Zhang \cite{Rapoport2017} 
(see \cite[\S14.1-14.2]{LZ} for the construction and more precise technical assumptions)
, which is semistable at all places of $E$ above $\mathcal{V}_\mathrm{ram}\cup \mathcal{V}_\mathrm{asd}$. When $K_G$ is the stabilizer of a global self-dual lattice, the regular integral model $\mathcal{X}_K$ recovers that in  \cite{Bruinier2017} if $F_0=\mathbb{Q}$. Let $\varphiK\in \sS(\mathbb{V}^n_f)^K$ be a factorizable Schwartz function such that $\varphi_{v}=\mathbf{1}_{(\Lambda_v)^n}$ at all $v$ nonsplit in $F$. Let $T\in \Herm_n(F_0)$ be nonsingular. Associated to $(T,\varphiK)$ we have an arithmetic special cycle $\mathcal{Z}(T,\varphiK)_K$ over $\XK$ (\cite[\S14.3]{LZ}). The definition of the arithmetic special cycle was pioneered by Kudla--Rapoport in both global and local settings (\cite{Kudla2014,Kudla2011}), and plays a fundamental role in the formulation of the arithmetic Siegel--Weil formula. 

\subsection{Modularity in arithmetic Chow groups}

The integral model $\XK$ and $\mathcal{Z}(T,\varphiK)_K$ are constructed as the moduli spaces of certain abelian varieties with additional structures. To describe them more precisely, in this subsection we consider the special case $F_0=\mathbb{Q}$ (so $E=F$ is an imaginary quadratic field), $\m=1$ and there is a global self-dual hermitian lattice $\Lambda$ such that $\Lambda_v=\Lambda \otimes_{O_{F_0}}O_{F_0,v}$ and $\varphi_v=\mathbf{1}_{\Lambda_v}$ at all finite places $v$. In this special case, the special cycles are indexed by $T\in\Herm_\m(O_{F_0})_{\ge0}=\mathbb{Z}_{\ge0}$.  Assume that $F/F_0$ is unramified at 2 for simplicity. 

Define an integral model $\XK$ of $\Sh_{K}$ over $O_F$ as follows. For an $O_F$-scheme $S$, we consider $\XK(S)$ to be the groupoid of tuples $(A_0,\iota_0,\lambda_0,A, \iota,\lambda,\mathcal{F}_A)$, where
\begin{enumerate}
\item $A_0$ (resp. $A$) is an abelian scheme over $S$.
\item $\iota_0$ (resp. $\iota$) is an action of $O_F$ on $A_0$ (resp. $A$).

\item $\lambda_0$ (resp. $\lambda$) is a principal polarization of $A_0$ (resp. $A$).
\item $\mathcal{F}_A\subseteq\Lie A$ is an $O_F$-stable $\mathcal{O}_S$-module local direct summand of rank $\n-1$.
\end{enumerate}
We require that
\begin{enumerate}
\item $O_F$ acts on the $\mathcal{O}_S$-module $\Lie A_0$ via the structure morphism $O_F\hookrightarrow \mathcal{O}_S$. This is the \emph{Kottwitz condition} of signature $(1,0)$: $$\det(T-\iota_0(a)|\Lie A_0)=T-a\in \mathcal{O}_S[T]$$ for $a\in O_F$.
\item $\mathcal{F}_A$ satisfies the \emph{Kr\"amer condition}: $O_F$ acts on $\mathcal{F}_A$ via the structure morphism and acts on the line bundle $\Lie A/\mathcal{F}_A$ via the conjugate of the structure morphism.   This implies (and in characteristic 0 is equivalent to) the Kottwitz condition of signature $(\n-1,1)$:  $$\det(T-\iota(a)|\Lie A)=(T-a)^{\n-1}(T-\bar a))\in \mathcal{O}_S[T]$$ for $a\in O_F$.
\item The Rosati involution on $\End A_0$ (resp. $\End A$) induces the conjugation on $O_F$ via $\iota_0$ (resp. $\iota$).
\item At every geometric point $s$ of $S$, there is an isomorphism of hermitian $O_{F,\ell}$-modules $$\Hom_{O_F}(T_\ell A_{0,s},T_\ell A_s)\simeq \Hom_{O_F}(\Lambda_0,\Lambda) \otimes \mathbb{Z}_\ell$$ for any prime $\ell$ different from the residue characteristic of $s$. Here $\Lambda_0$ is a fixed self-dual hermitian lattice of rank 1 over $O_F$. Notice that $\Hom_{O_F}(\Lambda_0,\Lambda)$ has a natural hermitian module structure given by $(x,y):=y^\vee\circ x\in\End_{O_F}(\Lambda_0)\subseteq F$ and similarly for the left-hand-side.
\end{enumerate}
Then the functor $S\mapsto \XK(S)$ is represented by a Deligne--Mumford stack $\XK$ regular over $\Spec O_{F}$.

The extra data $(A_0,\iota_0,\lambda_0)$ of a CM elliptic curve allows us to consider a motivic version of the lattice $\Lambda$. For $(A_0,\iota_0,\lambda_0,A, \iota,\lambda,\mathcal{F}_A)\in\XK(S)$, define the module of \emph{special homomorphisms} to be $$\Lambda(A_0,A):=\Hom_{O_F}(A_0,A),$$ equipped with a natural hermitian form $(x,y)\in O_F$ is given by $$(A_0\xrightarrow{x}A\xrightarrow{\lambda}A^\vee \xrightarrow{y^{\vee}}A_0^\vee\xrightarrow{\lambda_0^{-1}}A_0)\in\End_{O_F}(A_0)=\iota_0(O_F)\simeq O_F.$$

When $T>0$, define the special divisor $\mathcal{Z}(T,\varphiK)_K$ by requiring an additional special homomorphism of norm $T$. More precisely, the functor $S\mapsto \{(A_0,\iota_0,\lambda_0,A, \iota,\lambda,\mathcal{F}_A,x)\}$, where
\begin{enumerate}
\item $(A_0,\iota_0,\lambda_0,A, \iota,\lambda,\mathcal{F}_A)\in \XK(S)$,
\item $x\in \Lambda(A_0,A)$ such that $(x,x)=T$,
\end{enumerate}
is represented by a Deligne--Mumford stack  $\mathcal{Z}(T,\varphiK)_K$, which is finite and unramified over $\XK$.  It extends to a compactified special divisor $\mathcal{Z}^*(T,\varphiK)_K$ on the canonical toroidal compactification $\XK^*$ by taking the Zariski closure. \cite{Bruinier2017} further defines a \emph{total special divisor} $\mathcal{Z}^\mathrm{tot}(T,\varphiK)_K$ by adding an explicit boundary divisor to $\mathcal{Z}^*(T,\varphiK)_K$ (\cite[(1.1.3)]{Bruinier2017}). Using regularized theta lifts of harmonic Maass forms, $\mathcal{Z}^\mathrm{tot}(T,\varphiK)_K$ is equipped with an \emph{automorphic Green function} with log-log singularities along the boundary ((\cite[\S7.2]{Bruinier2017})), hence defines an element in $$\wh{\mathcal{Z}}^\mathrm{tot}(T,\varphiK)_K\in\aCh^1(\XK^*).$$

When $T=0$, define $$\wh{\mathcal{Z}}^\mathrm{tot}(0,\varphiK)_K=\wh{\mathcal{L}}_K^\vee+(\mathrm{Exc},-\log|\disc(F)|
)\in\aCh^1(\XK^*)$$ where $\wh{\mathcal{L}}_K^\vee$ is the metrized dual tautological line bundle over $\XK^*$, $\mathrm{Exc}$ is an effective vertical divisor supported above $\mathcal{V}_\mathrm{ram}$ equipped with the constant Green function $-\log|\disc(F)|
$. Define \emph{the generating function in arithmetic Chow groups} $$\wh{\mathcal{Z}}^\mathrm{tot}(\sz,\varphiK)_K:=\sum_{T\ge0}\wh{\mathcal{Z}}^\mathrm{tot}(T,\varphiK)_K\cdot q^T,$$ as a formal generating function valued in $\aCh^1(\XK^*)$, where $\sz\in \mathcal{H}_1$ lies in the usual upper half plane.

\begin{theorem}[Modularity in arithmetic Chow groups: the divisor case]\label{thm:BHKRY}
  The formal generating function $\wh{\mathcal{Z}}^\mathrm{tot}(\sz,\varphiK)_K$ 
  defines an elliptic modular form  valued in $\aCh^1(\XK^*)$ of weight $\n$, level $|\disc F|$ and character $\eta^\n$.  
\end{theorem}

This is proved in Bruinier--Howard--Kudla--Rapoport--Yang \cite[Theorem B]{Bruinier2017}. Analogous to \S\ref{sec:arithm-modul-conj}, Theorem \ref{thm:BHKRY} allows us to construct arithmetic theta lifts valued in $\aCh^1(\XK^*)$. As applications, \cite[Theorems A,B]{BHK+20} prove formulas relating the arithmetic intersection of these arithmetic theta lifts and small/big CM points to the central derivative of certain convolution $L$-functions of two elliptic modular forms, generalizing the Gross--Zagier formula \cite{Gross1986}.

\begin{remark}\label{rem:nonholo}
One can also use \emph{Kudla's Green function} \cite[(7.4.1)]{Bruinier2017} in place of the automorphic Green function to define an arithmetic divisor $\wh{\mathcal{Z}}^\mathrm{tot}(\sy, T,\varphiK)_K\in\aCh^1(\XK^*)$ depending on a parameter $\sy=\Im(\sz)\in \mathbb{R}_{>0}$ (here $T$ is also allowed to be $<0$, in which case the divisor is supported at the archimedean fiber). Then the generating function $$\sum_{T\in \mathbb{Z}} \wh{\mathcal{Z}}^\mathrm{tot}(\sy, T,\varphiK)_K\cdot q^T$$ becomes a nonholomorphic modular form (\cite[Theorem 7.4.1]{Bruinier2017}. This is a consequence of Theorem \ref{thm:BHKRY} and the modularity of the difference of the two generating functions due to Ehlen--Sankaran \cite{ES18}.
\end{remark}

\begin{remark}
The proof of Theorem \ref{thm:BHKRY} uses the arithmetic theory of Borcherds products, which requires the assumption $F_0=\mathbb{Q}$. For $F_0\ne \mathbb{Q}$, a version of Theorem \ref{thm:BHKRY} is proved in Qiu \cite{Qiu22} by a different method and formulation. A version of Theorem \ref{thm:BHKRY} is proved in Howard--Madapusi Pera \cite{HMP20} for (open) orthogonal Shimura varieties over $\mathbb{Q}$.  
\end{remark}

\begin{remark}
  The generating functions of arithmetic divisors have also found many applications outside the Kudla program. To name some recent arithmetic applications:
  \begin{enumerate}
  \item Theorem \ref{thm:BHKRY} is  used in Zhang's proof of the arithmetic fundamental lemma over $\mathbb{Q}_p$ in \cite{Zhang2019}. Variants over general totally real fields also play a key role for the arithmetic fundamental lemma over $p$-adic fields in Mihatsch--Zhang \cite{MZ21} and the arithmetic transfer conjecture in Z. Zhang \cite{Zha21}, within the framework of the arithmetic Gan--Gross--Prasad conjectures for unitary groups. We refer to Zhang's article in these proceedings for more details.
  \item The arithmetic modularity in \cite{HMP20} is used in Shankar--Shankar--Tang--Tayou \cite{SSTT2022} on the Picard rank jumps of K3 surfaces over number fields.
  \item The proof of the averaged Colmez conjecture in Andreatta--Goren--Howard--Madapusi Pera \cite{Andreatta2018} relies on relating arithmetic intersection of special divisors on orthogonal Shimura varities and big CM points to central derivatives of certain $L$-functions.
  \end{enumerate}
\end{remark}

\begin{remark}\label{rem:highcodim}
 The modularity problem in arithmetic Chow groups  \cite[Problem 4]{Kudla2004} remains open in higher codimension $\m>1$. When $\m>1$, even when $T>0$ the special cycle $\mathcal{Z}(T,\varphi)$ in general has the wrong codimension due to improper intersection in positive characteristics, and the consideration of \emph{derived intersection} is necessary to obtain the correct class $\wh{\mathcal{Z}}(T,\varphi)$ in arithmetic Chow groups. It is also subtle to find the correction terms at places of bad reduction and at boundary (both issues already appear when $\m=1$) and  to find the correct construction of Green currents to ensure modularity.

 The recent works of Howard--Madapusi \cite{HM22} and Madapusi \cite{Mad22} address some of these issues when $\m>1$. \cite{HM22} defines the special cycles of any codimension on the integral model of (open) Shimura varieties for orthogonal groups over $\mathbb{Q}$ and proves the modularity of their generating series. \cite{Mad22} uses methods from \emph{derived algebraic geometry} to define special cycles on the integral model of more general Hodge type Shimura varieties with good reduction, which recovers those constructed  using more classical methods in \cite{HM22}.
\end{remark}

\subsection{Arithmetic Siegel--Weil formula}

If the arithmetic theta function $\hZ\in\aCHXK$ can be constructed, then we may apply the  \emph{arithmetic volume} $$\avol: \aCHXK\rightarrow \mathbb{C}, \quad\wh{\mathcal{Z}}\mapsto\wdeg(\wh{\mathcal{Z}}\cdot (c_1(\wh{\mathcal{L}}_K^\vee))^{\dim \XK-\mm })$$ and try to relate $\avol(\hZ)$ to the special derivatives of Siegel Eisenstein series. However, as discussed in Remark \ref{rem:highcodim} the definition of $\hZ$ is rather subtle when $\m>1$. Moreover, the special derivatives are nonholomorphic modular forms, and thus for comparison it is better to construct nonholomorphic generating function (also including terms indexed by $T\not\in \Herm_\m(F_0)_{\ge0}$, cf. Remark \ref{rem:nonholo}).

In this subsection we assume that $m=n$, so $s_0=0$ and $\kappa=1$ in the Siegel--Weil formula for the pair $(V,W)$. In this special case, the arithmetic volume is simply the arithmetic degree and we can define the nonsingular terms in the generating function 
in a more explicit way. Even for $T>0$ terms, the relation to Siegel Eisenstein series is more complicated due to contribution at places of bad reduction, a phenomenon first discovered by Kudla--Rapoport  \cite{Kudla2000}   via explicit computation in the context of Shimura curves uniformized by the Drinfeld $p$-adic half plane.

 For nonzero $t_1,\ldots,t_n\in F_0$ and $\varphi_1,\cdots,\varphi_n\in\sS(\VAf)^K$ such  that $\varphi_{v}=\mathbf{1}_{\Lambda_v}$ at all $v$ nonsplit in $F$, we have a natural decomposition (cf. \cite[(11.2)]{Kudla2014}),
\begin{equation}
  \label{eq:Z(T)decomp}
  \mathcal{Z}(t_1,\varphi_1)_K\cap\cdots\cap\mathcal{Z}(t_n,\varphi_n))K=\bigsqcup_{T\in \Herm_{n}(F_0)} \mathcal{Z}(T, \varphiK)_K,
\end{equation}
 here $\cap$ denotes taking fiber product over $\XK$, the indexes $T$ have diagonal entries $t_1,\ldots,t_n$, and $\varphiK=\otimes_{i=1}^n \varphi_i$. For $v\nmid\infty$ and $\nu$ a place of $E$ over $v$,  define the \emph{local arithmetic intersection number}
\begin{equation}\label{eq:localint}
\Int_{T,\nu}(\varphiK)\coloneqq \chi(\mathcal{Z}(T,\varphiK)_K, \mathcal{O}_{\mathcal{Z}(t_1,\varphi_1)_K} \otimes^\mathbb{L}\cdots \otimes^\mathbb{L}\mathcal{O}_{\mathcal{Z}(t_n,\varphi_n)_K})\cdot\log q_\nu,  
\end{equation}
where $q_\nu$ denotes the size of the residue field $k_\nu$ of $E_{\nu}$, $\mathcal{O}_{\mathcal{Z}(t_i,\varphi_i)_K}$ denotes the structure sheaf of the special divisor $(\mathcal{Z}(t_i,\varphi_i)_K)_{O_{E_\nu}}$, $\otimes^\mathbb{L}$ denotes the derived tensor product of coherent sheaves on $(\XK)_{O_{E_\nu}}$, and $\chi$ denotes the Euler--Poincar\'e characteristic (an alternating sum of lengths of $\mathcal{O}_{E_\nu}$-modules). Notice that  the derived tensor product  $\mathcal{O}_{\mathcal{Z}(t_1,\varphi_1)_K} \otimes^\mathbb{L}\cdots \otimes^\mathbb{L}\mathcal{O}_{\mathcal{Z}(t_n,\varphi_n)_K}$ has the structure of a complex of $\mathcal{O}_{\mathcal{Z}(t_1,\varphi_1)_K\cap\cdots \cap\mathcal{Z}(t_n,\phi_n)_K}$-modules, hence has a natural decomposition by support according to the decomposition (\ref{eq:Z(T)decomp}).  Define $$\Int_{T,v}(\varphiK)\coloneqq \frac{1}{[E:F_0]}\cdot \sum_{\nu|v}\Int_{T,\nu}(\varphiK).$$  Using the star product of \emph{Kudla's Green functions}, we can also define its local arithmetic intersection number  $\Int_{T,v}(\sy,\varphiK)$ at infinite places (\cite[\S15.3]{LZ}), which depends on a parameter $\sy\in \Herm_n(F_{0,\infty})_{>0}$. Combining all the local arithmetic numbers together, define the \emph{global arithmetic intersection number}, or the (normalized) \emph{arithmetic degree} of the special cycle $\mathcal{Z}(T,\varphiK)_K$ 
\begin{equation}
  \label{eq:wchdecomp}
  \wdeg_T(\sy,\varphiK)\coloneqq\frac{1}{\vol([\Sh_K])/2}\left( \sum_{v\nmid\infty}\Int_{T,v}(\varphiK)+\sum_{v\mid \infty}\Int_{T,v}(\sy,\varphiK)\right).
\end{equation}

We form the \emph{generating function of arithmetic degrees}  $$\wdeg(\sz, \varphiK)\coloneqq \sum_{T\in\Herm_n(F_0)\atop \det T\ne0}\wdeg_T(\sy,\varphiK) q^T.$$

On the other hand, associated to $$\varphiI:=\varphiK \otimes \varphi_\infty\in\sS(\mathbb{V}^n),$$ where $\varphi_{\infty}$ is the Gaussian function, we obtain a classical \emph{incoherent Eisenstein series} $E(\sz, s,\varphiI)$ (\cite[\S12.4]{LZ}, this terminology originates from Kudla \cite[\S2]{Kudla1997a}). The central value $E(\sz, 0, \varphiI)=0$ by the incoherence. We thus consider its \emph{central derivative} $$\Eis'(\sz, \varphi)\coloneqq \frac{\rd}{\rd s}\bigg|_{s=0}E(\sz, s,\varphiI).$$

To match the arithmetic degree, we need to modify $\Eis'(\sz, \varphiK)$ by central values of coherent Eisenstein series at places of bad reduction. For $v\in\mathcal{V}_\mathrm{ram}\cup\mathcal{V}_\mathrm{asd}$, let  $^v\mathbb{V}$ be the \emph{coherent} hermitian space over $\mathbb{A}_F$ nearby $\mathbb{V}$ at $v$, namely $(^v\mathbb{V})_w\simeq \mathbb{V}_w$ exactly for all places $w\ne v$. For any vertex lattice $\Lambda_{t,v}\subseteq (^v\mathbb{V})_v$ of type $t$, the Schwartz function $\varphi^v \otimes \mathbf{1}_{(\Lambda_{t,v})^n}\otimes \varphi_\infty\in\sS((^v\mathbb{V})^n)$ gives a classical \emph{coherent Eisenstein series} $E(\sz,s,\varphi^v \otimes \mathbf{1}_{(\Lambda_{t,v})^n}\otimes \varphi_\infty)$. Define the (normalized) \emph{central values} $$^v\Eis_t(\sz,\varphiK):=\frac{\vol(K_{G,v})}{\vol(K_{\Lambda_{t,v}})}\cdot E(\sz,0,\varphi^v \otimes \mathbf{1}_{(\Lambda_{t,v})^n}\otimes \varphi_\infty).$$ 
Define the \emph{modified central derivative} $$\pEis(\sz,\varphiK):=\Eis'(\sz,\varphiK)+(-1)^n\sum_{v\in\mathcal{V}_\mathrm{ram}\cup\mathcal{V}_\mathrm{asd}}{}^v\Eis(\sz,\varphiK).$$ Here ${}^v\Eis(\sz,\varphiK)$ is an explicit $\mathbb{Q}$-linear combination of $^v\Eis_t(\sz,\varphiK)$ for certain $t$'s as defined in \cite{HLSY23}. 
It has a decomposition into Fourier coefficients $$\pEis(\sz,\varphiK)=\sum_{T\in \Herm_n(F_0)}\pEis_T(\sz,\varphiK).$$

Now we can state the arithmetic Siegel--Weil formula, which is an identity between the arithmetic degrees and 
the modified central derivative of the incoherent Eisenstein series.

\begin{theorem}[Arithmetic Siegel--Weil formula: nonsingular terms] \label{thm:ASW}  Assume that $F/F_0$ is split at all places above 2.  Let $\varphiK\in \sS(\mathbb{V}^n_f)^K$ be a factorizable Schwartz function such that $\varphi_{v}=\mathbf{1}_{(\Lambda_v)^n}$ at all $v$ nonsplit in $F$. Let $T\in \Herm_n(F_0)$ be nonsingular. Then
$$\wdeg_T(\sy, \varphiK)q^T=(-1)^n\cdot \pEis_T(\sz,\varphiK).$$ Further assume that $\varphiK$ is nonsingular (\cite[\S12.3]{LZ}) at two places split in $F$. Then $$\wdeg(\sz, \varphiK)=(-1)^n\cdot \pEis(\sz,\varphiK).$$ In particular, $\wdeg(\sz, \varphiK)$ is a nonholomorphic hermitian modular form on $\mathcal{H}_n$.
\end{theorem}

The proof of this theorem boils down to a local arithmetic Siegel--Weil formula computing $\Int_{T,v}(\varphi)$ at each place $v$ nonsplit in $F$:
\begin{enumerate}
\item At $v\mid\infty$, this is the archimedean arithmetic Siegel--Weil formula proved by  Liu \cite{Liu2011} and Garcia--Sankaran \cite{Garcia2019} independently.
\item At $v\nmid\infty$ inert in $F$ such that $\Lambda_v$ is self-dual, this is the content of the \emph{Kudla--Rapoport conjecture}  ({\cite[Conjecture 11.10]{Kudla2014}}), proved by Zhang and the author \cite{LZ}. We refer to \cite[\S5]{Li23} for an exposition. An analogous theorem is also proved for orthogonal Shimura varieties over $\mathbb{Q}$ at a place of good reduction \cite{LZ22}.
\item At $v\nmid\infty$ inert in $F$ such that $\Lambda_v$ is almost self-dual, this is a variant of the Kudla--Rapoport conjecture formulated and proved by Zhang and the author \cite{LZ}.
\item At $v\nmid\infty$ ramified in $F$ such that $\Lambda_v$ is self-dual, this is the Kudla--Rapoport conjecture for Kr\"amer models formulated by He--Shi--Yang  \cite{HSY2023} and proved by He--Shi--Yang and the author \cite{HLSY23}.    
\end{enumerate}

\begin{remark}
 The precise formulation of the singular part of the arithmetic Siegel--Weil \cite[Problem 6]{Kudla2004} remains an open problem. As a special case, the constant term of the arithmetic Siegel--Weil formula should roughly relate the arithmetic volume of $\XK$ to logarithmic derivatives of Dirichlet $L$-functions. Such an explicit arithmetic volume formula is proved by Bruinier--Howard \cite{BH21}, though a precise comparison with the constant term of $\pEis(\sz,\varphiK)$ is yet to be formulated and established.
 \end{remark}

 \begin{remark}
   In contrast to the classical and geometric Siegel--Weil formula, in Theorem \ref{thm:ASW} the choice of $K\subseteq H(\mathbb{A}_f)$ is fixed at all nonsplit places $v$ in order to construct a regular integral model $\XK$, which prevents us from formulating a full adelic version of the arithmetic Siegel--Weil formula. Nevertheless, the flexibility at split places (due to the regular integral models with Drinfeld level structure) allow us to choose $\varphi$ to be nonsingular at split places to kill all singular terms on both sides. This extra flexibility is crucial for applications such as the arithmetic inner product formula, via making use of the multiplicity one result of the doubling method and bypassing the need for proving the singular part of the arithmetic Siegel--Weil formula.

   A related open problem is to formulate and prove an arithmetic Siegel--Weil formula when the level $K$ is more general at nonsplit places. For minuscule parahoric levels at inert places, such a formulation can be found in Cho \cite{Cho2022}.
  \end{remark}

  \begin{remark}
  Over function fields,  the recent work of Feng--Yun--Zhang \cite{FYZ2024} proved a higher Siegel--Weil formula  for unitary groups in the unramified setting, which relates nonsingular coefficients of the $r$-th derivative of Siegel Eisenstein series and intersection numbers of special cycles on moduli spaces of Drinfeld shtukas with $r$ legs.  The case $r=0$ (resp. $r=1$) can be viewed as an analogue of the Siegel--Weil formula (resp. the arithmetic Siegel--Weil formula).  Over number fields, however, no analogue of such a higher Siegel--Weil formula is currently known when $r>1$. Feng--Yun--Zhang  \cite{FYZ21a} further defined higher theta series over function fields (including all singular terms), using both classical and derived algebraic geometry, and conjectured their modularity (see Feng--Harris's article \cite[\S4]{FH2023} in these proceedings for more details). 
  \end{remark}

\subsection{Arithmetic inner product formula} In this subsection we come back to the equal rank situation consider in \S\ref{sec:theta-dich-equal} and assume that $m=2n$, so we can take $\chi$ to be the trivial character, and $s_0=0$, $\kappa=1$ in the Siegel--Weil formula for the pair $(V,W^\square)$. Assume $\varepsilon(\pi)=-1$ so $\mathbb{V}=\mathbb{V}_\pi$ is incoherent. In this case we may use the Beilinson--Bloch height pairing to define an arithmetic inner product between arithmetic theta lifts.

To this end, assume that $V$ is anisotropic, thus $X_K$ is projective. Let  $\Ch^\mm(X_K)^0\subseteq \Ch^\mm(X_K)$ be the subgroup of cohomologically trivial cycles.  Since $\dim X_K=2n-1$ we have a (conditional) symmetric bilinear height pairing   $$\langle\ ,\  \rangle_\mathrm{BB}: \Ch^\mm(X_K)^0\times \Ch^{\mm}(X_K)^0\rightarrow \mathbb{R},$$ constructed Beilinson \cite{Beuilinson1987} and Bloch \cite{Bloch1984a}. It generalizes the Neron--Tate height pairing when $\mm=1$. When $\m>1$, the Beilinson--Bloch height pairing is only defined assuming certain conjectures on algebraic cycles on $X_K$ (see \cite[Conjectures 2.2.1 and 2.2.3]{Beuilinson1987}). This important technical issue is addressed in \cite{LL2020, LL2021} so that the left-hand-side of (\ref{eq:AIPF}) in Theorem \ref{thm:AIPF} can be defined unconditionally, but we will intentionally ignore it for the purpose of this article. Then we naturally obtain an inner product on $\CHXK^0$ and define the (normalized) \emph{arithmetic inner product} $$\langle\ ,\  \rangle_{X_K}: \CHXK^0\times\CHXK^0\rightarrow \mathbb{C},\quad(Z_1,Z_2)\mapsto\frac{\langle Z_1,Z_2\rangle_\mathrm{BB}}{\vol([X_K])/2},$$ which also gives a well-defined inner product $\langle\ ,\ \rangle_X$ on $\CHX^0$.

  \begin{assumption}\label{ass:pi} We impose the following (mild) local assumptions on $F/F_0$ and $\pi$.
    \begin{enumerate}
      \item $F/F_0$ is split at all 2-adic places and $F_0\ne \mathbb{Q}$. If $v\nmid \infty$ is ramified in $F$, then $v$ is unramified over $\mathbb{Q}$. Assume that $F/\mathbb{Q}$ is Galois or contains an imaginary quadratic field (for simplicity).
      \item For $v|\infty$, $\pi_v$ is the holomorphic discrete series with Harish-Chandra parameter $\{\frac{\n-1}{2}, \frac{\n-3}{2},\allowbreak\ldots,\frac{-\n+3}{2}, \frac{-\n+1}{2}\}$.
      \item For $v\nmid\infty$, $\pi_v$ is tempered.
      \item For $v\nmid\infty$ ramified in $F$, $\pi_v$ is spherical with respect to the stabilizer of $O_{F_v}^{2\mm}$.
      \item\label{item:inert} For $v\nmid\infty$ inert in $F$, $\pi_v$ is unramified or almost unramified (\cite{Liu22}) with respect to the stabilizer of $O_{F_v}^{2\mm}$. If $\pi_v$ is almost unramified, then $v$ is unramified over $\mathbb{Q}$.
      \end{enumerate}
    \end{assumption}

    Under Assumption \ref{ass:pi}, the arithmetic theta lift $\Theta_\varphi(\phi)$ is in fact cohomologically trivial and thus $\Theta_\mathbb{V}(\pi)\subseteq\CHX^0$ (see \cite[Proposition 6.10]{LL2020}) and we can apply the arithmetic inner product.

          \begin{theorem}[Arithmetic inner product formula]\label{thm:AIPF}
            Let $\pi$ be a cuspidal automorphic representation of $G(\mathbb{A})$ satisfying Assumption \ref{ass:pi}. Assume that $\varepsilon(\pi)=-1$. Assume that Kudla's modularity Conjecture~\ref{conj:modularity} holds. Then for any $\phi_i=\otimes_v\phi_{i,v}\in \pi\cap\sA_{\n}(G(\mathbb{A}))$, $\varphi_i=\otimes_v\varphi_{i,v}\in \sS(\VAf^{n})$ ($i=1,2$), 
            \begin{equation}
              \label{eq:AIPF}
               \langle\Theta_{\varphi_1}(\phi_1),\Theta_{\varphi_2}(\phi_2)\rangle_{X}=\frac{L'(1/2,\pi)}{b_{2n}(0)}\cdot \prod_v Z_v^\natural(0,\phi_{1,v},\phi_{2,v}, \varphiI_{1,v},\varphiI_{2,v}).
            \end{equation} In particular, $$L'(1/2,\pi)\ne0\Longrightarrow \Theta_{\mathbb{V}}(\pi)\ne0,$$ and the converse also holds if $\langle\ ,\ \rangle_X$ is nondegenerate.
          \end{theorem}

This is proved by Liu and the author in \cite{LL2020, LL2021}. The conjectural arithmetic inner product formula was formulated (in the orthogonal case) by Kudla \cite{Kudla1997a} using the Gillet--Soul\'e height and in more generality by Liu \cite{Liu2011} using the Beilinson--Bloch height. This theorem verifies (under local assumptions) the conjecture formulated by Liu (who also completely proved the case $\m=1$ in \cite{Liu2011a}). 

          \begin{remark}
            The formula can further be made explicit by computing the local doubling zeta integrals. For example, if
    \begin{itemize}
    \item $\pi$ is unramified or almost unramified at all finite places,
    \item $\phi\in\pi$ is a holomorphic newform such that $(\phi,\overline{\phi})_\pi=1$,
    \item $\varphi_v$ is the characteristic function of self-dual or almost self-dual lattices at all finite places $v$.
    \end{itemize}
    Then we have
    \begin{equation}\label{eq:explicit}
\langle \Theta_\varphi(\phi),\Theta_\varphi(\phi)\rangle_{X}=(-1)^\m\frac{L'(1/2,\pi)}{b_{2\m}(0)}C_\m^{[F_0:\mathbb{Q}]}\prod_{v\in S_\pi}\frac{q_v^{\m-1}(q_v+1)}{(q_v^{2\m-1}+1)(q_v^{2\m}-1)},
\end{equation}
where $C_\m=2^{-2\m}\pi^{\m^2}\frac{\Gamma(1)\cdots\Gamma(\m)}{\Gamma(\m+1)\cdots\Gamma(2\m)}$ is an archimedean doubling zeta integral computed by Eischen--Liu \cite{EL20} and  $S_\pi=\{v \text{ inert}: \pi_v\text{ almost unramified}\}$.

Notice that the Grand Riemann Hypothesis predicts that $L'(1/2,\pi)\ge0$, while Beilinson's Hodge index conjecture (\cite[Conjecture 5.5]{Beuilinson1987}) predicts that $(-1)^\m\langle \Theta_\varphi(\phi),\Theta_\varphi(\phi)\rangle_X\ge0$. It is a good reality check that these two (big) conjectures are compatible with \eqref{eq:explicit}.         \end{remark}

The arithmetic inner product formula can be viewed as a higher dimensional generalization of the Gross--Zagier formula \cite{Gross1986}, and has applications to the Beilinson--Bloch conjecture for higher dimensional Shimura varieties. Without assuming Kudla's modularity conjecture, we cannot define $\Theta_\varphi(\phi)$ but we may still obtain unconditional nonvanishing results  (\cite{LL2020, LL2021}) on Chow groups as predicted by the Beilinson--Bloch conjecture (using a proof by contradiction argument).

\begin{theorem}[Application to the Beilinson--Bloch conjecture]\label{thm:BB}
  Let $\pi$ be a cuspidal automorphic representation of $G(\mathbb{A})$ satisfying Assumption \ref{ass:pi}. Let $\Ch^{\m}(X)^0_{\mathfrak{m}_\pi}$ the localization of $\Ch^{\m}(X)^0_\mathbb{C}$ at the maximal ideal $\mathfrak{m}_\pi$ of the spherical Hecke algebra of $H(\mathbb{A}_f)$ (away from all ramification) associated to $\pi$. Then the implication
  \begin{equation}
    \label{eq:BBmpi}
    \ord_{s=1/2}L(s,\pi)=1\Longrightarrow \rank\Ch^\m(X)^0_{\mathfrak{m}_\pi}\ge1
  \end{equation}
  holds when the level subgroup $K\subseteq H(\mathbb{A}_f)$ is sufficiently small.
\end{theorem}
    
\begin{remark}
  Disegni--Liu \cite{DL22} proved a $p$-adic version of the arithmetic inner product formula, relating the central derivative of the cyclotomic $p$-adic $L$-function $L_p(\pi)$ to the $p$-adic height pairing (\cite{Nek93}) of arithmetic theta lifts. As an application, they proved implications of the form
  \begin{center}
  central order of vanishing of $L_p(\pi)$ is 1 $\Longrightarrow$ Bloch--Kato Selmer group $\H^1_f(F,\rho_\pi(\m))$ has rank $\ge1$,
  \end{center} where $\rho_\pi$ is the Galois representation associated to $\pi$. This verifies part of the $p$-adic Bloch--Kato conjecture.
\end{remark}

\begin{remark}
  Xue \cite{Xue19} used the arithmetic inner product formula in the case $\m=1$ and the Gan--Gross--Prasad conjecture for $\U(2)\times\U(2)$ to prove endoscopic cases of the arithmetic Gan--Gross--Prasad conjecture for $\U(2)\times \U(3)$. In general, one also expects a similar relation between the arithmetic inner product formula for $\U(\n)$ and endoscopic cases of arithmetic Gan--Gross--Prasad conjecture for $\U(\n)\times\U(\n+1)$.
\end{remark}

\begin{remark}
  Throughout we have assumed the skew-hermitian space $W$ has even dimension $2n$. When the skew-hermitian space $W$ has odd dimension, Liu \cite{Liu21b} has defined mixed arithmetic theta lifting and used it to formulate a conjectural arithmetic inner product formula.
\end{remark}

\subsection{Summary}We end with a summary of the trilogy of theta correspondences in Table~\ref{tab:summary}.

\begin{table}[h]
  \centering\setlength\extrarowheight{2pt}
  \begin{tabular}{c|c|c|c}
      & Theta function & Siegel--Weil formula & Inner product formula\\\hline 
      Classical  & $\theta(g,h,\varphi)$ & $I(g,\varphi)\stackrel{\cdot}{=} E(g,s_0,\varphi)$ & $\langle\theta_{\varphi}(\phi),\theta_{\varphi}(\phi)\rangle_{H}\ceq L(s_0+\frac{1}{2},\pi)$ \\
& (\S\ref{sec:theta-lifting}) & (Theorem \ref{thm:SWF}) & (Theorem \ref{thm:rallis-inner-product}) \\\hline
      Geometric & $[Z(g,\varphi)]$ & $\voln[Z(g,\varphi)]\ceq E(g,s_0,\varphic)$ & $\langle\theta^\KM_{\varphi}(\phi),\theta^\KM_{\varphi}(\phi)\rangle_{X(\mathbb{C})}\ceq L(s_0+\frac{1}{2},\pi)$ \\
& (\S\ref{sec:geometric-modularity}) & (Theorem \ref{thm:geoSWF}) & (Theorem \ref{thm:geom-inner-prod}) \\\hline
      Arithmetic  & $Z(g,\varphi)$ & $\wdeg(\sz, \varphiK)\ceq \pEis(\sz,\varphi)$ & $\langle\Theta_{\varphi}(\phi),\Theta_{\varphi}(\phi)\rangle_{X}\ceq L'(\frac{1}{2},\pi)$\\
& (\S\ref{sec:kudl-gener-funct}) & (Theorem \ref{thm:ASW}) & (Theorem \ref{thm:AIPF})
    \end{tabular}
  \caption{Summary}
  \label{tab:summary}
\end{table}

\bibliographystyle{alpha}
\bibliography{KR}

\newcommand{\etalchar}[1]{$^{#1}$}
\begin{thebibliography}{AGHMP18}

\bibitem[AGHMP18]{Andreatta2018}
Fabrizio Andreatta, Eyal~Z. Goren, Benjamin Howard, and Keerthi Madapusi~Pera.
\newblock Faltings heights of abelian varieties with complex multiplication.
\newblock {\em Ann. of Math. (2)}, 187(2):391--531, 2018.

\bibitem[Bei87]{Beuilinson1987}
A.~A. Beilinson.
\newblock Height pairing between algebraic cycles.
\newblock In {\em {$K$}-theory, arithmetic and geometry ({M}oscow,
  1984--1986)}, volume 1289 of {\em Lecture Notes in Math.}, pages 1--25.
  Springer, Berlin, 1987.

\bibitem[BF10]{BF10}
Jan~Hendrik Bruinier and Jens Funke.
\newblock On the injectivity of the {K}udla-{M}illson lift and surjectivity of
  the {B}orcherds lift.
\newblock In {\em Moonshine: the first quarter century and beyond}, volume 372
  of {\em London Math. Soc. Lecture Note Ser.}, pages 12--39. Cambridge Univ.
  Press, Cambridge, 2010.

\bibitem[BGKK07]{BGKK07}
J.~I. Burgos~Gil, J.~Kramer, and U.~K\"{u}hn.
\newblock Cohomological arithmetic {C}how rings.
\newblock {\em J. Inst. Math. Jussieu}, 6(1):1--172, 2007.

\bibitem[BH21]{BH21}
Jan~Hendrik {Bruinier} and Benjamin {Howard}.
\newblock {Arithmetic volumes of unitary Shimura varieties}.
\newblock {\em arXiv e-prints}, page arXiv:2105.11274, May 2021.

\bibitem[BHK{\etalchar{+}}20a]{Bruinier2017}
Jan~H. Bruinier, Benjamin Howard, Stephen~S. Kudla, Michael Rapoport, and
  Tonghai Yang.
\newblock Modularity of generating series of divisors on unitary {S}himura
  varieties.
\newblock {\em Ast\'{e}risque}, (421, Diviseurs arithm\'{e}tiques sur les
  vari\'{e}t\'{e}s orthogonales et unitaires de Shimura):7--125, 2020.

\bibitem[BHK{\etalchar{+}}20b]{BHK+20}
Jan~H. Bruinier, Benjamin Howard, Stephen~S. Kudla, Michael Rapoport, and
  Tonghai Yang.
\newblock Modularity of generating series of divisors on unitary {S}himura
  varieties {II}: {A}rithmetic applications.
\newblock {\em Ast\'{e}risque}, (421, Diviseurs arithm\'{e}tiques sur les
  vari\'{e}t\'{e}s orthogonales et unitaires de Shimura):127--186, 2020.

\bibitem[Blo84]{Bloch1984a}
Spencer Bloch.
\newblock Height pairings for algebraic cycles.
\newblock In {\em Proceedings of the {L}uminy conference on algebraic
  {$K$}-theory ({L}uminy, 1983)}, volume~34, pages 119--145, 1984.

\bibitem[BMM16]{BMM16}
Nicolas Bergeron, John Millson, and Colette Moeglin.
\newblock The {H}odge conjecture and arithmetic quotients of complex balls.
\newblock {\em Acta Math.}, 216(1):1--125, 2016.

\bibitem[Bor99]{Bor99}
Richard~E. Borcherds.
\newblock The {G}ross-{K}ohnen-{Z}agier theorem in higher dimensions.
\newblock {\em Duke Math. J.}, 97(2):219--233, 1999.

\bibitem[Bru12]{Bru12}
Jan~Hendrik Bruinier.
\newblock Regularized theta lifts for orthogonal groups over totally real
  fields.
\newblock {\em J. Reine Angew. Math.}, 672:177--222, 2012.

\bibitem[BWR15]{BW15}
Jan~Hendrik Bruinier and Martin Westerholt-Raum.
\newblock Kudla's modularity conjecture and formal {F}ourier-{J}acobi series.
\newblock {\em Forum Math. Pi}, 3:e7, 30, 2015.

\bibitem[BZ22]{BZ22}
Jan~Hendrik Bruinier and Shaul Zemel.
\newblock Special cycles on toroidal compactifications of orthogonal {S}himura
  varieties.
\newblock {\em Math. Ann.}, 384(1-2):309--371, 2022.

\bibitem[Cho22]{Cho2022}
Sungyoon Cho.
\newblock Special cycles on unitary {S}himura varieties with minuscule
  parahoric level structure.
\newblock {\em Math. Ann.}, 384(3-4):1747--1813, 2022.

\bibitem[DL22]{DL22}
Daniel {Disegni} and Yifeng {Liu}.
\newblock {A $p$-adic arithmetic inner product formula}.
\newblock {\em arXiv e-prints}, page arXiv:2204.09239, April 2022.

\bibitem[Dun22]{Dun22}
Nguyen~Chi Dung.
\newblock {\em Geometric {P}ullback {F}ormula for {U}nitary {S}himura
  {V}arieties}.
\newblock ProQuest LLC, Ann Arbor, MI, 2022.
\newblock Thesis (Ph.D.)--Columbia University.

\bibitem[EGT23]{EGT23}
Philip {Engel}, Fran{\c{c}}ois {Greer}, and Salim {Tayou}.
\newblock {Mixed mock modularity of special divisors}.
\newblock {\em arXiv e-prints}, page arXiv:2301.05982, January 2023.

\bibitem[EL20]{EL20}
Ellen {Eischen} and Zheng {Liu}.
\newblock {Archimedean Zeta Integrals for Unitary Groups}.
\newblock {\em arXiv e-prints}, page arXiv:2006.04302, June 2020.

\bibitem[ES18]{ES18}
Stephan Ehlen and Siddarth Sankaran.
\newblock On two arithmetic theta lifts.
\newblock {\em Compos. Math.}, 154(10):2090--2149, 2018.

\bibitem[FH23]{FH2023}
Tony Feng and Michael Harris.
\newblock {Derived structures in the Langlands correspondence}.
\newblock 2023.
\newblock \url{https://math.berkeley.edu/~fengt/Feng-Harris.pdf}.

\bibitem[FM06]{FM06}
Jens Funke and John Millson.
\newblock Cycles with local coefficients for orthogonal groups and
  vector-valued {S}iegel modular forms.
\newblock {\em Amer. J. Math.}, 128(4):899--948, 2006.

\bibitem[FM14]{FM14}
Jens Funke and John Millson.
\newblock The geometric theta correspondence for {H}ilbert modular surfaces.
\newblock {\em Duke Math. J.}, 163(1):65--116, 2014.

\bibitem[FYZ21]{FYZ21a}
Tony {Feng}, Zhiwei {Yun}, and Wei {Zhang}.
\newblock {Higher theta series for unitary groups over function fields}.
\newblock {\em arXiv e-prints}, page arXiv:2110.07001, October 2021.

\bibitem[FYZ24]{FYZ2024}
Tony Feng, Zhiwei Yun, and Wei Zhang.
\newblock Higher {S}iegel--{W}eil formula for unitary groups: the non-singular
  terms.
\newblock {\em Invent. Math.}, 235(2):569--668, 2024.

\bibitem[Gan23]{Gan2023}
Wee~Teck Gan.
\newblock {Explicit Constructions of Automorphic Forms: Theta Correspondence
  and Automorphic Descent}.
\newblock {\em arXiv e-prints}, page arXiv:2303.14919, March 2023.

\bibitem[GI14]{GI14}
Wee~Teck Gan and Atsushi Ichino.
\newblock Formal degrees and local theta correspondence.
\newblock {\em Invent. Math.}, 195(3):509--672, 2014.

\bibitem[GKZ87]{GKZ87}
B.~Gross, W.~Kohnen, and D.~Zagier.
\newblock Heegner points and derivatives of {$L$}-series. {II}.
\newblock {\em Math. Ann.}, 278(1-4):497--562, 1987.

\bibitem[GQT14]{GQT14}
Wee~Teck Gan, Yannan Qiu, and Shuichiro Takeda.
\newblock The regularized {S}iegel-{W}eil formula (the second term identity)
  and the {R}allis inner product formula.
\newblock {\em Invent. Math.}, 198(3):739--831, 2014.

\bibitem[Gro21]{Gro21}
Benedict~H. Gross.
\newblock Incoherent definite spaces and shimura varieties.
\newblock In Werner M{\"u}ller, Sug~Woo Shin, and Nicolas Templier, editors,
  {\em Relative Trace Formulas}, pages 187--215, Cham, 2021. Springer
  International Publishing.

\bibitem[GS90]{GS90}
Henri Gillet and Christophe Soul\'{e}.
\newblock Arithmetic intersection theory.
\newblock {\em Inst. Hautes \'{E}tudes Sci. Publ. Math.}, (72):93--174 (1991),
  1990.

\bibitem[GS19]{Garcia2019}
Luis~E. Garcia and Siddarth Sankaran.
\newblock Green forms and the arithmetic {S}iegel-{W}eil formula.
\newblock {\em Invent. Math.}, 215(3):863--975, 2019.

\bibitem[GZ86]{Gross1986}
Benedict~H. Gross and Don~B. Zagier.
\newblock Heegner points and derivatives of {$L$}-series.
\newblock {\em Invent. Math.}, 84(2):225--320, 1986.

\bibitem[HK91]{HK1991}
Michael Harris and Stephen~S. Kudla.
\newblock The central critical value of a triple product {$L$}-function.
\newblock {\em Ann. of Math. (2)}, 133(3):605--672, 1991.

\bibitem[HKS96]{Harris1996}
Michael Harris, Stephen~S. Kudla, and William~J. Sweet.
\newblock Theta dichotomy for unitary groups.
\newblock {\em J. Amer. Math. Soc.}, 9(4):941--1004, 1996.

\bibitem[HLSY23]{HLSY23}
Qiao He, Chao Li, Yousheng Shi, and Tonghai Yang.
\newblock A proof of the {K}udla-{R}apoport conjecture for {K}r\"{a}mer models.
\newblock {\em Invent. Math.}, 234(2):721--817, 2023.

\bibitem[HM22]{HM22}
Benjamin {Howard} and Keerthi {Madapusi}.
\newblock {Kudla's modularity conjecture on integral models of orthogonal
  Shimura varieties}.
\newblock {\em arXiv e-prints}, page arXiv:2211.05108, November 2022.

\bibitem[HMP20]{HMP20}
Benjamin Howard and Keerthi Madapusi~Pera.
\newblock Arithmetic of {B}orcherds products.
\newblock {\em Ast\'{e}risque}, (421, Diviseurs arithm\'{e}tiques sur les
  vari\'{e}t\'{e}s orthogonales et unitaires de Shimura):187--297, 2020.

\bibitem[HSY23]{HSY2023}
Qiao He, Yousheng Shi, and Tonghai Yang.
\newblock Kudla-{R}apoport conjecture for {K}r\"{a}mer models.
\newblock {\em Compos. Math.}, 159(8):1673--1740, 2023.

\bibitem[HZ76]{Hirzebruch1976}
F.~Hirzebruch and D.~Zagier.
\newblock Intersection numbers of curves on {H}ilbert modular surfaces and
  modular forms of {N}ebentypus.
\newblock {\em Invent. Math.}, 36:57--113, 1976.

\bibitem[Ich07]{Ich07}
Atsushi Ichino.
\newblock On the {S}iegel-{W}eil formula for unitary groups.
\newblock {\em Math. Z.}, 255(4):721--729, 2007.

\bibitem[KM86]{Kudla1986}
Stephen~S. Kudla and John~J. Millson.
\newblock The theta correspondence and harmonic forms. {I}.
\newblock {\em Math. Ann.}, 274(3):353--378, 1986.

\bibitem[KM87]{KM87}
Stephen~S. Kudla and John~J. Millson.
\newblock The theta correspondence and harmonic forms. {II}.
\newblock {\em Math. Ann.}, 277(2):267--314, 1987.

\bibitem[KM90]{KM90}
Stephen~S. Kudla and John~J. Millson.
\newblock Intersection numbers of cycles on locally symmetric spaces and
  {F}ourier coefficients of holomorphic modular forms in several complex
  variables.
\newblock {\em Inst. Hautes \'{E}tudes Sci. Publ. Math.}, (71):121--172, 1990.

\bibitem[KMSW14]{KMSW14}
Tasho {Kaletha}, Alberto {Minguez}, Sug~Woo {Shin}, and Paul-James {White}.
\newblock {Endoscopic Classification of Representations: Inner Forms of Unitary
  Groups}.
\newblock {\em arXiv e-prints}, page arXiv:1409.3731, September 2014.

\bibitem[KR88]{KR88}
Stephen~S. Kudla and Stephen Rallis.
\newblock On the {W}eil-{S}iegel formula.
\newblock {\em J. Reine Angew. Math.}, 387:1--68, 1988.

\bibitem[KR94]{KR94}
Stephen~S. Kudla and Stephen Rallis.
\newblock A regularized {S}iegel-{W}eil formula: the first term identity.
\newblock {\em Ann. of Math. (2)}, 140(1):1--80, 1994.

\bibitem[KR00]{Kudla2000}
Stephen~S. Kudla and Michael Rapoport.
\newblock Height pairings on {S}himura curves and {$p$}-adic uniformization.
\newblock {\em Invent. Math.}, 142(1):153--223, 2000.

\bibitem[KR11]{Kudla2011}
Stephen~S. Kudla and Michael Rapoport.
\newblock Special cycles on unitary {S}himura varieties {I}. {U}nramified local
  theory.
\newblock {\em Invent. Math.}, 184(3):629--682, 2011.

\bibitem[KR14]{Kudla2014}
Stephen~S. Kudla and Michael Rapoport.
\newblock Special cycles on unitary {S}himura varieties {II}: {G}lobal theory.
\newblock {\em J. Reine Angew. Math.}, 697:91--157, 2014.

\bibitem[KRY06]{Kudla2006}
Stephen~S. Kudla, Michael Rapoport, and Tonghai Yang.
\newblock {\em Modular forms and special cycles on {S}himura curves}, volume
  161 of {\em Annals of Mathematics Studies}.
\newblock Princeton University Press, Princeton, NJ, 2006.

\bibitem[KSZ]{Kisina}
Mark Kisin, Sug~Woo Shin, and Yihang Zhu.
\newblock Cohomology of certain {Shimura} varieties of abelian type (temporary
  title).
\newblock in preparation.

\bibitem[Kud94]{Kudla1994}
Stephen~S. Kudla.
\newblock Splitting metaplectic covers of dual reductive pairs.
\newblock {\em Israel J. Math.}, 87(1-3):361--401, 1994.

\bibitem[Kud97a]{Kudla1997}
Stephen~S. Kudla.
\newblock Algebraic cycles on {S}himura varieties of orthogonal type.
\newblock {\em Duke Math. J.}, 86(1):39--78, 1997.

\bibitem[Kud97b]{Kudla1997a}
Stephen~S. Kudla.
\newblock Central derivatives of {E}isenstein series and height pairings.
\newblock {\em Ann. of Math. (2)}, 146(3):545--646, 1997.

\bibitem[Kud03]{Kud03}
Stephen~S. Kudla.
\newblock Integrals of {B}orcherds forms.
\newblock {\em Compositio Math.}, 137(3):293--349, 2003.

\bibitem[Kud04]{Kudla2004}
Stephen~S. Kudla.
\newblock Special cycles and derivatives of {E}isenstein series.
\newblock In {\em Heegner points and {R}ankin {$L$}-series}, volume~49 of {\em
  Math. Sci. Res. Inst. Publ.}, pages 243--270. Cambridge Univ. Press,
  Cambridge, 2004.

\bibitem[Kud21]{Kud2021}
Stephen~S. Kudla.
\newblock Remarks on generating series for special cycles on orthogonal
  {S}himura varieties.
\newblock {\em Algebra Number Theory}, 15(10):2403--2447, 2021.

\bibitem[Lap08]{Lap2008}
Erez~M. Lapid.
\newblock A remark on {E}isenstein series.
\newblock In {\em Eisenstein series and applications}, volume 258 of {\em
  Progr. Math.}, pages 239--249. Birkh\"{a}user Boston, Boston, MA, 2008.

\bibitem[{Lap}22]{Lap2022}
Erez {Lapid}.
\newblock {Some perspectives on Eisenstein series}.
\newblock {\em arXiv e-prints}, page arXiv:2204.02914, April 2022.

\bibitem[Li92]{Li92}
Jian-Shu Li.
\newblock Nonvanishing theorems for the cohomology of certain arithmetic
  quotients.
\newblock {\em J. Reine Angew. Math.}, 428:177--217, 1992.

\bibitem[Li23]{Li23}
Chao Li.
\newblock From sum of two squares to arithmetic {S}iegel-{W}eil formulas.
\newblock {\em Bull. Amer. Math. Soc. (N.S.)}, 60(3):327--370, 2023.

\bibitem[Liu11a]{Liu2011}
Yifeng Liu.
\newblock Arithmetic theta lifting and {$L$}-derivatives for unitary groups,
  {I}.
\newblock {\em Algebra Number Theory}, 5(7):849--921, 2011.

\bibitem[Liu11b]{Liu2011a}
Yifeng Liu.
\newblock Arithmetic theta lifting and {$L$}-derivatives for unitary groups,
  {II}.
\newblock {\em Algebra Number Theory}, 5(7):923--1000, 2011.

\bibitem[Liu21]{Liu21b}
Yifeng Liu.
\newblock Mixed arithmetic theta lifting for unitary groups.
\newblock In Werner M{\"u}ller, Sug~Woo Shin, and Nicolas Templier, editors,
  {\em Relative Trace Formulas}, pages 329--350, Cham, 2021. Springer
  International Publishing.

\bibitem[Liu22]{Liu22}
Yifeng Liu.
\newblock Theta correspondence for almost unramified representations of unitary
  groups.
\newblock {\em J. Number Theory}, 230:196--224, 2022.

\bibitem[LL21]{LL2020}
Chao Li and Yifeng Liu.
\newblock Chow groups and {$L$}-derivatives of automorphic motives for unitary
  groups.
\newblock {\em Ann. of Math. (2)}, 194(3):817--901, 2021.

\bibitem[LL22]{LL2021}
Chao Li and Yifeng Liu.
\newblock Chow groups and {$L$}-derivatives of automorphic motives for unitary
  groups, {II}.
\newblock {\em Forum of Mathematics, Pi}, 10:e5, 71pp, 2022.

\bibitem[LR05]{LR05}
Erez~M. Lapid and Stephen Rallis.
\newblock On the local factors of representations of classical groups.
\newblock In {\em Automorphic representations, {$L$}-functions and
  applications: progress and prospects}, volume~11 of {\em Ohio State Univ.
  Math. Res. Inst. Publ.}, pages 309--359. de Gruyter, Berlin, 2005.

\bibitem[LZ22a]{LZ}
Chao Li and Wei Zhang.
\newblock Kudla-{R}apoport cycles and derivatives of local densities.
\newblock {\em J. Amer. Math. Soc.}, 35(3):705--797, 2022.

\bibitem[LZ22b]{LZ22}
Chao Li and Wei Zhang.
\newblock On the arithmetic {S}iegel-{W}eil formula for {GS}pin {S}himura
  varieties.
\newblock {\em Invent. Math.}, 228(3):1353--1460, 2022.

\bibitem[{Mad}22]{Mad22}
Keerthi {Madapusi}.
\newblock {Derived special cycles on Shimura varieties}.
\newblock {\em arXiv e-prints}, page arXiv:2212.12849, December 2022.

\bibitem[Mae21]{Mae2021}
Yota Maeda.
\newblock The modularity of special cycles on orthogonal {S}himura varieties
  over totally real fields under the {B}eilinson-{B}loch conjecture.
\newblock {\em Canad. Math. Bull.}, 64(1):39--53, 2021.

\bibitem[Mae22]{Mae2022}
Yota Maeda.
\newblock Modularity of special cycles on unitary {S}himura varieties over
  {CM}-fields.
\newblock {\em Acta Arith.}, 204(1):1--18, 2022.

\bibitem[Mok15]{Mok15}
Chung~Pang Mok.
\newblock Endoscopic classification of representations of quasi-split unitary
  groups.
\newblock {\em Mem. Amer. Math. Soc.}, 235(1108):vi+248, 2015.

\bibitem[{Mor}23]{Morel2023}
Sophie {Morel}.
\newblock {Shimura varieties}.
\newblock {\em arXiv e-prints}, page arXiv:2310.16184, October 2023.

\bibitem[MZ21]{MZ21}
Andreas {Mihatsch} and Wei {Zhang}.
\newblock {On the Arithmetic Fundamental Lemma conjecture over a general
  $p$-adic field}.
\newblock {\em arXiv e-prints}, page arXiv:2104.02779, April 2021.

\bibitem[Nek93]{Nek93}
Jan Nekov\'{a}\v{r}.
\newblock On {$p$}-adic height pairings.
\newblock In {\em S\'{e}minaire de {T}h\'{e}orie des {N}ombres, {P}aris,
  1990--91}, volume 108 of {\em Progr. Math.}, pages 127--202. Birkh\"{a}user
  Boston, Boston, MA, 1993.

\bibitem[Pra90]{Pra1990}
Dipendra Prasad.
\newblock Trilinear forms for representations of {${\rm GL}(2)$} and local
  {$\epsilon$}-factors.
\newblock {\em Compositio Math.}, 75(1):1--46, 1990.

\bibitem[PSR86]{PR86}
I.~Piatetski-Shapiro and S.~Rallis.
\newblock {$\epsilon$} factor of representations of classical groups.
\newblock {\em Proc. Nat. Acad. Sci. U.S.A.}, 83(13):4589--4593, 1986.

\bibitem[PSR87]{PR87}
I.~Piatetski-Shapiro and S.~Rallis.
\newblock {\em L-functions for the classical groups}, pages 1--52.
\newblock Springer Berlin Heidelberg, Berlin, Heidelberg, 1987.

\bibitem[{Qiu}22]{Qiu22}
Congling {Qiu}.
\newblock {Arithmetic modularity of special divisors and arithmetic mixed
  Siegel-Weil formula (with an appendix by Yujie Xu)}.
\newblock {\em arXiv e-prints}, page arXiv:2204.13457, April 2022.

\bibitem[Ral84]{Ral84}
S.~Rallis.
\newblock Injectivity properties of liftings associated to {W}eil
  representations.
\newblock {\em Compositio Math.}, 52(2):139--169, 1984.

\bibitem[Ral87]{Rallis1987}
Stephen Rallis.
\newblock {\em {$L$}-functions and the oscillator representation}, volume 1245
  of {\em Lecture Notes in Mathematics}.
\newblock Springer-Verlag, Berlin, 1987.

\bibitem[RSZ20]{Rapoport2017}
M.~Rapoport, B.~Smithling, and W.~Zhang.
\newblock Arithmetic diagonal cycles on unitary {S}himura varieties.
\newblock {\em Compos. Math.}, 156(9):1745--1824, 2020.

\bibitem[Sou92]{Sou92}
C.~Soul\'{e}.
\newblock {\em Lectures on {A}rakelov geometry}, volume~33 of {\em Cambridge
  Studies in Advanced Mathematics}.
\newblock Cambridge University Press, Cambridge, 1992.
\newblock With the collaboration of D. Abramovich, J.-F. Burnol and J. Kramer.

\bibitem[SSTT22]{SSTT2022}
Ananth~N. Shankar, Arul Shankar, Yunqing Tang, and Salim Tayou.
\newblock Exceptional jumps of {P}icard ranks of reductions of {K}3 surfaces
  over number fields.
\newblock {\em Forum Math. Pi}, 10:Paper No. e21, 49, 2022.

\bibitem[Tun83]{Tun1983}
Jerrold~B. Tunnell.
\newblock Local {$\epsilon $}-factors and characters of {${\rm GL}(2)$}.
\newblock {\em Amer. J. Math.}, 105(6):1277--1307, 1983.

\bibitem[Wal85]{Wal1985}
J.-L. Waldspurger.
\newblock Sur les valeurs de certaines fonctions {$L$} automorphes en leur
  centre de sym\'{e}trie.
\newblock {\em Compositio Math.}, 54(2):173--242, 1985.

\bibitem[Wei65]{Weil1965}
Andr\'{e} Weil.
\newblock Sur la formule de {S}iegel dans la th\'{e}orie des groupes
  classiques.
\newblock {\em Acta Math.}, 113:1--87, 1965.

\bibitem[Xia22]{Xia22}
Jiacheng Xia.
\newblock Some cases of {K}udla's modularity conjecture for unitary {S}himura
  varieties.
\newblock {\em Forum Math. Sigma}, 10:Paper No. e37, 31, 2022.

\bibitem[Xue19]{Xue19}
Hang Xue.
\newblock Arithmetic theta lifts and the arithmetic {G}an-{G}ross-{P}rasad
  conjecture for unitary groups.
\newblock {\em Duke Math. J.}, 168(1):127--185, 2019.

\bibitem[Yam11]{Yam11}
Shunsuke Yamana.
\newblock On the {S}iegel-{W}eil formula: the case of singular forms.
\newblock {\em Compos. Math.}, 147(4):1003--1021, 2011.

\bibitem[Yam14]{Yam14}
Shunsuke Yamana.
\newblock L-functions and theta correspondence for classical groups.
\newblock {\em Invent. Math.}, 196(3):651--732, 2014.

\bibitem[YZZ09]{YZZ09}
Xinyi Yuan, Shou-Wu Zhang, and Wei Zhang.
\newblock The {G}ross-{K}ohnen-{Z}agier theorem over totally real fields.
\newblock {\em Compos. Math.}, 145(5):1147--1162, 2009.

\bibitem[Zha09]{Zha09}
Wei Zhang.
\newblock {\em Modularity of generating functions of special cycles on
  {S}himura varieties}.
\newblock ProQuest LLC, Ann Arbor, MI, 2009.
\newblock Thesis (Ph.D.)--Columbia University.

\bibitem[Zha19]{Zha19}
Shouwu Zhang.
\newblock Linear forms, algebraic cycles, and derivatives of {L}-series.
\newblock {\em Sci. China Math.}, 62(11):2401--2408, 2019.

\bibitem[Zha21a]{Zhang2019}
Wei Zhang.
\newblock Weil representation and {A}rithmetic {F}undamental {L}emma.
\newblock {\em Ann. of Math. (2)}, 193(3):863--978, 2021.

\bibitem[{Zha}21b]{Zha21}
Zhiyu {Zhang}.
\newblock {Maximal parahoric arithmetic transfers, resolutions and modularity}.
\newblock {\em arXiv e-prints}, page arXiv:2112.11994, December 2021.

\end{thebibliography}

\end{document}